\documentclass[UTF-8,reqno]{amsart}
\usepackage{enumerate, bbm}
\usepackage{amssymb,url,color, booktabs}
\usepackage{mathrsfs}
\usepackage{soul}
\usepackage[nobysame,abbrev]{amsrefs}
\BibSpec{article}{%
+{}{\PrintAuthors} {author}
+{,}{ \textrm} {title}
+{.}{ \textit} {journal}
+{,}{ \textbf} {volume}
+{}{ \parenthesize} {date}
+{,}{ } {pages}
+{.}{ To appear in \textit} {status}
+{.}{ arXiv:} {eprint}
+{.}{} {transition}
}
\BibSpec{book}{%
+{}{\PrintAuthors} {author}
+{,}{ \textit} {title}
+{.}{ \textrm} {series} 
+{,}{ Vol.} {volume} 
+{.}{ } {publisher}
+{,}{ } {date}
+{,}{ } {pages}
+{.}{} {transition}
}
\usepackage{color}
\usepackage[colorlinks=true]{hyperref}
\hypersetup{
    linkcolor=blue,          
    citecolor=red,        
    filecolor=blue,      
    urlcolor=cyan
}

\usepackage{color}
\definecolor{MyDarkBlue}{cmyk}{0.8,0.3,0.8,0.4}
\definecolor{yellow}{rgb}{0.99,0.99,0.70}
\definecolor{white}{rgb}{1.0,1.0,1.0}
\definecolor{black}{rgb}{0.00,0.00,0.00}
\definecolor{backgroundcolor}{RGB}{199,238,206}

\numberwithin{equation}{section}

\newcommand{\be}{\begin{eqnarray}}
\newcommand{\ee}{\end{eqnarray}}
\newcommand{\ce}{\begin{eqnarray*}}
\newcommand{\de}{\end{eqnarray*}}
\newtheorem{theorem}{Theorem}[section]
\newtheorem{lemma}[theorem]{Lemma}
\newtheorem{remark}[theorem]{Remark}
\newtheorem{definition}[theorem]{Definition}
\newtheorem{proposition}[theorem]{Proposition}
\newtheorem{Examples}[theorem]{Example}
\newtheorem{corollary}[theorem]{Corollary}

\def\tr{\mathrm {tr}}

\def\e{\mathrm{e}}
\def\supp{\mathrm{supp}}

\def\dif{{\mathord{{\rm d}}}}

\def\bba{{\boldsymbol{a}}}

\def\bbk{{\boldsymbol{k}}}

\def\bb2{{\boldsymbol{2}}}
\def\no{\nonumber}
\def\={&\!\!=\!\!&}

\def\bx{{\mathbf{x}}}

\def\bB{{\mathbf B}}
\def\bC{{\mathbf C}}
\def\bE{{\mathbf E}}

\def\bP{{\mathbf P}}

\def\b1{{\mathbbm 1}}

\def\cM{{\mathcal M}}

\def\cR{{\mathcal R}}

\def\mB{{\mathbb B}}
\def\mC{{\mathbb C}}

\def\mE{{\mathbb E}}

\def\mI{{\mathbb I}}

\def\mH{{\mathbb H}}
\def\mL{{\mathbb L}}

\def\mN{{\mathbb N}}

\def\mP{{\mathbb P}}
\def\mQ{{\mathbb Q}}
\def\mR{{\mathbb R}}

\def\sB{{\mathscr B}}

\def\sF{{\mathscr F}}

\def\sI{{\mathscr I}}

\def\sP{{\mathscr P}}

\def\sS{{\mathscr S}}

\def\geq{\geqslant}
\def\leq{\leqslant}
\def\<{{\langle}}
\def\>{{\rangle}}
\def\({{\big(}}
\def\){{\big)}}
\def\[{{\Big[}}
\def\]{{\Big]}}

\def\a{\alpha}
\def\b{\beta}
\def\de{\delta}
\def\g{\gamma}
\def\G{\Gamma}
\def\k{\kappa}
\def\l{\lambda}
\def\om{\omega}

\def\s{\sigma}

\def\t{\theta}

\def\ff{\frac}
\def\nn{\nabla}
\def\ss{\sqrt}
\def\p{\partial}
\def\div{\mathord{{\rm div}}}

\def\bt{\begin{theorem}}
\def\et{\end{theorem}}
\def\bl{\begin{lemma}}
\def\el{\end{lemma}}
\def\br{\begin{remark}}
\def\er{\end{remark}}
\def\bpf{\begin{proof}}
\def\epf{\end{proof}}
\def\bx{\begin{Examples}}
\def\ex{\end{Examples}}
\def\bd{\begin{definition}}
\def\ed{\end{definition}}
\def\bp{\begin{proposition}}
\def\ep{\end{proposition}}
\def\bc{\begin{corollary}}
\def\ec{\end{corollary}}

\def\wt{\widetilde}

\allowdisplaybreaks

\begin{document}

\title{Kinetic SDEs with subcritical distributional drifts}
\author{Zikai Chen, Zimo Hao and Xicheng Zhang}


\thanks{
This work is supported by National Key R\&D program of China (No. 2023YFA1010103) and NNSFC grant of China (No. 12131019)  and the DFG through the CRC 1283 ``Taming uncertainty and profiting from randomness and low regularity in analysis, stochastics and their applications''. Z. Chen is also supported by the China Scholarship Council (Grant 202506270058).
}

\begin{abstract}
In this paper we study the well-posedness of the kinetic stochastic differential equation (SDE) in $\mR^{2d}(d\geq2)$ driven by Brownian motion:
$$\dif X_t=V_t\dif t,\ \dif V_t=b(t,X_t,V_t)\dif t+\sqrt{2}\dif W_t,$$
where the subcritical distribution-valued drift $b$ belongs to the weighted anisotropic H\"{o}lder space $\mL_T^{q_b}\bC_\bba^{\a_b}(\rho_\k)$ with parameters $\a_b\in(-1,0)$, $q_b\in(\ff{2}{1+\a_b},\infty]$, $\k\in[0,1+\a_b)$ and $\div_v b$ is bounded. We establish the well-posedness of weak solutions to the associated integral equation:
$$X_t=X_0+\int_0^t V_s\dif s,\ V_t=V_0+\lim_{n\to\infty}\int_0^t b_n(s,X_s,V_s)\dif s+\sqrt{2}W_t,$$
where $b_n:=b*\G_n$ denotes the mollification of $b$ and the limit is taken in the $L^2$-sense. 
As an application, we discuss examples of $b$ involving Gaussian random fields.
\end{abstract}

\maketitle

\section{Introduction}
Fix $d\geq2$ and $T>0$. In this paper, we investigate the following second-order stochastic differential equation driven by standard Brownian motion $W_t$ in $\mR^{d}$:
\begin{align}\label{1}
\ddot{X_t}=b(t,X_t,\dot{X_t})+\sqrt{2}\dot{W_t},\ t\in[0,T],
\end{align}
where the drift $b$ is a time-dependent distribution. A typical example is Gaussian random fields. In this case, $X_t$ models the evolution of particles in random environment (see \cite{Fan98,HZ23} for the first order case). By introducing the velocity variable $V_t=\dot{X_t}$, the SDE \eqref{1} can be written as the following first-order degenerate system: for $t\geq0$,
\begin{align}
\left\{
\begin{aligned}\label{2}
&\dif X_t=V_t\dif t,\\
&\dif V_t=b(t,X_t,V_t)\dif t+\sqrt{2}\dif W_t.
\end{aligned}
\right.
\end{align}
Let $Z_t:=(X_t,V_t)$. For $z=(x,v)$, we define 
$$
B(t,z)=B(t,x,v):=(v,b(t,x,v))^{\rm t}\in\mR^{2d},\ \sigma:=(0,\sqrt{2}\mI)^{\rm t}\in\mR^{2d}\otimes\mR^d,
$$ 
where ${\rm t}$ denotes the transpose of a matrix. Then the SDE \eqref{2} can be rewritten as
\begin{align}\label{3}
\dif Z_t=B(t,Z_t)\dif t+\sigma\dif  W_t,\ Z_0\sim\mu\in\sP(\mR^{2d}),
\end{align}
where $\sP(\mR^{2d})$ is the space of all probability measures on $\mR^{2d}$.

The main challenging part of our article is that the drift $b$ is distribution-valued; hence, the drift term is meaningless in the classical sense as we cannot assign a value to a distribution at the point $(X_t,V_t)$. To solve this problem, as in \cite{HZ23}, we employ a mollifying approximation to define a weak solution. Let $\Gamma_n(z)=\Gamma_n(x,v):=n^{4d}\Gamma(nx,n^3v)$, where $\Gamma\in C_c^\infty(\mR^{2d})$ is a smooth probability density function with compact support. We define the approximation of $b$ as follows:
\begin{align}
b_n(t,z):=(b(t,\cdot)*\Gamma_n)(z),\ B_n(t,z)=B_n(t,x,v):=(v,b_n(t,x,v))^{\rm t},
\end{align}
where $*$ denotes convolution in the distributional sense.

To begin, we introduce the concept of weak solutions to the above degenerate system with a distributional drift. 

\bd[Weak solutions]\label{weaksol}
Let $\mathfrak F:=(\Omega,\sF,(\sF_t)_{t\geq0},\bP)$ be a stochastic basis, and $(Z,W)$ be a pair of $\mR^{2d}$-valued continuous $\sF_t$-adapted processes on $\mathfrak F$. We call $(\mathfrak F,Z,W)$ a weak solution of SDE \eqref{3} with initial distribution $\mu\in\sP(\mR^{2d})$ if $W$ is an $\sF_t$-Brownian motion, $\bP\circ Z_0^{-1}=\mu$, and
\begin{align}\label{weaksolappro}
Z_t=Z_0+A_t^B+\s W_t,\ t\in[0,T],\ a.s.,
\end{align}
where $A_t^B:=\lim_{n\to\infty}\int_0^t B_n(s,Z_s)\dif s$ exists in the $L^2$-sense.
\ed

It should be noted that the definition of a weak solution in Definition \ref{weaksol} relies on the choice of mollifiers $\G_n$ (see Remark \ref{rk1.3} below). Our goal is to determine how weak the condition on $b$ can be for SDE \eqref{3} to have a solution in the sense of Definition \ref{weaksol}. 

\subsection{Main results}

To state our main results, we first introduce the following notations. Define the scaling parameter $\bba:=(3,1)$ and the corresponding anisotropic distance: for $z=(x,v)$ and $z'=(x',v')$ in $\mR^{2d}$,
\begin{align}\label{distance}
|z-z'|_\bba:=|x-x'|^{\ff{1}{3}}+|v-v'|.
\end{align}
We need the following weight function as well: for $\k\in\mR$,
\begin{align}\label{rhok}
    \rho_\k(z):=((1+|x|^2)^{1/3}+1+|v|^2)^{-\k/2}\asymp(1+|z|_\bba)^{-\k},\ z=(x,v)\in\mR^{2d}.
\end{align}
It is easy to see that for any $\kappa\in\mR$ and $j\in\mN$, there exists a constant $C=C(\k,j,d)>0$ such that
\begin{align}\label{1.14}
    |\nabla^j_v\rho_\kappa(z)|\lesssim_C\rho_\kappa(z)\rho_j(z),\ |\nabla^j_x\rho_\kappa(z)|\lesssim_C\rho_\kappa(z)\rho_{3j}(z).
\end{align}
For any $\kappa,\alpha\in\mR$, we define
\begin{align}\label{0117:01}
\|f\|_{\bC^\a_\bba(\rho_\k)}:=\|\rho_\kappa f\|_{\bC_\bba^\a},
\end{align}
where $\bC_\bba^\a$ denotes the anisotropic H\"{o}lder space (see Definition \ref{Holder} below for precise definition). The following theorems are our main results:

\bt\label{sec1:main1}
Let $\a_b\in(-1,0)$, $q_b\in(\ff{2}{1+\a_b},\infty]$ and $\k\in[0,\ff{1+\a_b-2/q_b}{3+\a_b-2/q_b}]$. Suppose that 
$$b\in\mL_T^{q_b}\bC^{\a_b}_\bba(\rho_\k),\ \div_v b\in\mL_T^{q_b}\bC^{\a_b}_{\bba}.$$
Let $\a\in(\ff{2}{q_b}+\ff{3\k-1}{1-\k},\a_b]$ and $q\in(\ff{2-2\k}{1+(1-\k)\a-3\k},q_b]$. Then for any $\de\in\mR$, $f\in\mL_T^q\bC_{\bba}^\a(\rho_\de)$ and $\l\geq0$, there exists a unique weak solution to PDE
\begin{align}\label{1:sec1:main1}
\p_t u=\Delta_v u-v\cdot\nn_x u-\l u+b\cdot\nn_v u+f,\ u_0\equiv0,
\end{align}
such that for any $\t\in[0,2-\ff{2}{q})$,
\begin{align}
(1+\l)^{\ff{\t}{2}}\|u\|_{\mL_T^\infty\bC_{\bba}^{2+\a-2/q-\t}(\rho_\nu)}\lesssim_C \|f\|_{\mL_T^q\bC_{\bba}^\a(\rho_\de)},\no
\end{align}
where $\nu:=\ff{2\k}{1+\a-2/q}+\de$ and the constant $C>0$ is independent of $\l$.
\et

\br
The detailed version of Theorem \ref{sec1:main1} is given in Theorem \ref{pdees}. 
We observe that when $\k=\de=0$, the estimate reduces to the classical Schauder's estimate (cf. \cite[Theorem 1.3]{CHM21}, \cite[Theorem 2.17]{HRZ23}). 
This theorem demonstrates that the growth or decay behavior of $u$ is jointly determined by both the growth of $b$ and the growth or decay of $f$.
\er

\bt\label{sec1:main2}
Let $\a_b\in(-1,0)$, $q_b\in(\ff{2}{1+\a_b},\infty]$ and $\k\in[0,1+\a_b)$. Suppose that $b\in\mL_T^{q_b}\bC^{\a_b}_\bba(\rho_\k)$ has a bounded divergence in $v$. Then for any $z_0\in\mR^{2d}$ and $p\geq2$,
there exists a unique weak solution $(\mathfrak F,Z,W)$ to SDE \eqref{3} starting from $z_0$ such that the following Krylov's estimate holds: for any $f\in\mL_T^1 C_b\cap\mL_T^q\bC_\bba^\a$ with $\a\in(-1,\a_b]$, $q\in(\ff{2}{1+\a},q_b]$ and $0\leq t_0<t_1\leq T$,
\begin{align}\label{sec1:main2:es1}
\Big\|\int_{t_0}^{t_1}f(s,Z_s)\dif s\Big\|_{L^p(\Omega)}\lesssim_C (t_1-t_0)^{\ff{2+\a-2/q}{2}}\|f\|_{\mL_T^q\bC_\bba^\a},
\end{align}
where the constant $C$ is independent of $t_0$, $t_1$ and $f$. Additionally, we have the following moment estimate: for any $\de\in\mR$,
\begin{align}
\bE\bigg(\sup_{t\in[0,T]}\rho_\de(Z_t)\bigg)\lesssim_{C(\de)}\rho_\de(z_0).\no
\end{align}
Furthermore, if $\k\in[0,\ff{1+\a_b-2/q_b}{3+\a_b-2/q_b}]$, then for any $f\in\mL_T^1 C_b(\rho_\de)\cap\mL_T^q\bC_\bba^\a(\rho_\de)$ with $\a\in(\ff{2}{q_b}+\ff{3\k-1}{1-\k},\a_b]$, $q\in(\ff{2-2\k}{1+(1-\k)\a-3\k},q_b]$ and $\de\in\mR$, we have for $\nu:=\ff{2\k}{1+\a-2/q}+\de$,
\begin{align}\label{sec1:main2:es2}
\Big\|\int_{t_0}^{t_1}f(s,Z_s)\dif s\Big\|_{L^p(\Omega)}\lesssim_C \rho_{-\nu}(z_0)(t_1-t_0)^{\ff{2+\a-2/q}{2}}\|f\|_{\mL_T^q\bC_{\bba}^\a(\rho_\de)},
\end{align}
where the constant $C$ is also independent of $t_0$, $t_1$ and $f$.
\et

\br\label{rk1.3}
We note that, as stated in Proposition \ref{limit} below, the solution mentioned above does not depend on the choice of mollifiers. Furthermore, for $A^B_t$ defined in \eqref{weaksolappro}, the function $t \mapsto A^B_t$ admits a H\"{o}lder continuous version with index $\beta<\ff{2+\a-2/q}{2}$.
\er

\br
The Krylov's estimate \eqref{sec1:main2:es1} implies that for any $z_0\in\mR^{2d}$ and Lebesgue almost all $t>0$, the process $Z_t(z_0)$ admits a density $\rho(t,z_0,z)$. Furthermore, by \eqref{sec1:main2:es2} we have
\begin{align}
\int_0^T\int_{\mR^{2d}}f(t,z)\rho(t,z_0,z)\dif z\dif t\lesssim\rho_{-\nu}(z_0)\|f\|_{\mL_T^q\bC_{\bba}^\a(\rho_\de)}.
\end{align}
Hence, by duality, $\rho$ belongs to the weighted space $\mL_T^{q'}\bB_{1,\bba}^{-\a,1}(\rho_{-\de})$ (see \cite[Proposition 2.76]{BCD11}), exhibiting certain growth or decay properties. 
\er

\br
In fact, by a similar argument, the well-posedness remains valid for any initial distribution $\mu$ with finite moment of order
$m$, provided that $$m(1+\a_b-\tfrac{2}{q_b})>2.$$
\er

We summarize our main contributions as follows:
\begin{itemize}
\item In Theorem \ref{sec1:main1} we establish the well-posedness of equation \eqref{1:sec1:main1} with the drift term $b$ belonging to the weighted anisotropic H\"{o}lder space $\mL_T^{q_b}\bC^{\a_b}_\bba(\rho_\k)$. This result extends the analysis to cases where the drift may be represented by certain Gaussian random fields.
\item While previous work \cite{HZZZ24} established weak well-posedness for $\alpha_b\in(-\frac{2}{3},-\frac{1}{2})$ with Gaussian noise $b$ (without requiring $\div_v b$ assumptions) by using paracontrolled calculus and renormalization techniques, Theorem \ref{sec1:main2} relaxes these conditions. Our result extends the admissible range to $\a_b\in(-1,0)$ and eliminates the probabilistic constraints on $b$ imposed in \cite[Theorem 6.3]{HZZZ24}, though it requires additional conditions on $\div_v b$.
\end{itemize}

\subsection{Background and motivation}

The kinetic equation, originally introduced by Landau in 1936, has been extensively studied in various fields such as physics, chemistry and biology, as it describes dynamic behaviors influenced by both deterministic laws and stochastic forces. For a general background introduction, we refer to \cite{Van Kampen11}; for ergodicity of the system, see \cite{HN04} and \cite{Villani09}; and for Euler approximations of invariant measures, see \cite{MSH02}.

For degenerate system with singular drifts:
\begin{align}
\label{nondedi}
\dif X_t=V_t\dif t,\ \ \dif V_t=b(t,X_t,V_t)\dif t+\sqrt{2}\sigma(t,X_t,V_t)\dif W_t,
\end{align}
many research has been conducted for the case where $b$ is a real function. In \cite{Chau17}, Chaudru de Raynel first studied equation \eqref{nondedi} under certain H\"{o}lder condition. In a subsequent work \cite{WZ16}, Wang and Zhang extended this result under a similar H\"{o}lder-Dini condition. In \cite{FFPV17}, the authors derived strong well-posedness for SDE \eqref{nondedi} with constant diffusion. Later, in \cite{Zhang18}, under weaker drift conditions than \cite{FFPV17}, Zhang established the strong well-posedness for SDE \eqref{nondedi} when $b$ belongs to anisotropic Sobolev spaces. In \cite{Zhang21}, Zhang derived the maximum principle for SDE \eqref{nondedi} with a distribution-valued inhomogeneous term. In a recent work \cite{CM22}, Chaudru de Raynel and Menozzi achieved the weak well-posedness under much weaker assumptions. In \cite{GRW24}, the authors established the exponential ergodicity of SDE \eqref{nondedi}. Additionally, kinetic SDEs driven by $\a$-stable process were studied in \cite{HWZ20}. However, it is important to note that all the aforementioned studies exclusively focus on drifts with positive regularity. For the case of negative regularity drifts, that is, when the drift $b$ is distribution-valued, defining solutions to the SDE \eqref{nondedi} presents significant challenges. In a recent work \cite{IPRT24}, the authors established the well-posedness of SDE \eqref{nondedi} for the regularity $\a\in(-\ff{1}{2},0)$ by using a classical method of the martingale problem.

Here we consider the kinetic system with degenerate diffusion \eqref{3} and categorize it into three cases through a simple scaling analysis. Let $\dot{\bC}_\bba^\a$ denote the homogeneous anisotropic H\"{o}lder space. Suppose that for some $q\in[1,\infty]$ and $\a\in\mR$, $b\in L^q(\mR_+; \dot{\bC}_\bba^\a)$, and that SDE \eqref{3} has a solution $Z$. For $\l>0$, we define 
$$Z_t^\l=(X_t^\l,V_t^\l):=(\l^{-3}X_{\l^2t},\l^{-1}V_{\l^2t}),\ W_t^\l:=\l^{-1}W_{\l^2t}.$$
Formally, it follows that 
$$\dif Z_t^\l=B^\l(t,Z_t^\l)\dif t+\s\dif W_t^\l,$$
where 
$$
B^\l(t,z)=B^\l(t,x,v):=(v,b^\l(t,x,v))^{\rm t}:=(v, \l b(\l^2 t,\l^3x,\l v))^{\rm t}.
$$ 
By the change of variables, we have
$$\|b^\l\|_{L^q(\mR_+;\dot{\bC}_\bba^\a)}=\l^{1+\a-\ff{2}{q}}\|b\|_{L^q(\mR_+; \dot{\bC}_\bba^\a)}.$$
Letting $\l\rightarrow0$, we shall categorize the following three cases:
$$\textbf{Subcritical:\ }\tfrac{2}{q}<1+\a;\textbf{\ Critical:\ }\tfrac{2}{q}=1+\a;\textbf{\ Supercritical:\ }\tfrac{2}{q}>1+\a.$$

Our aim is to establish the well-posedness of SDE \eqref{3} under the subcritical condition. To solve SDE \eqref{3}, we need a regularity estimate for the associated Kolmogorov equation
\begin{align}\label{Koleq}
\p_t u=\Delta_v u-v\cdot\nn_x u-\l u+b\cdot\nn_v u+f.
\end{align}
The study of kinetic Kolmogorov equations has been extensively explored in the literature. For example, in \cite{CHM21}, Chaudru de Raynel, Honor\'{e} and Menozzi derived sharp Schauder's estimates for equation \eqref{Koleq} under the condition that the coefficients belong to certain anisotropic H\"{o}lder spaces and that the first order term is nonlinear and unbounded. Here we focus on the case where $b\in\bC_\bba^\a$ with $\a<0$. According to the Schauder's theory of kinetic equation (see \cite{HRZ23}), the solution $u$ belongs to at most $\bC_\bba^{2+\a}$. In the classical sense, to ensure that the product $b\cdot\nn_v u$ is well-defined, we need $1+2\a>0$, which implies that $\a>-\ff{1}{2}$. This class of problems also arises in the study of singular SPDEs, such as in the analysis of the KPZ equation. To date, two principal methodologies have emerged as powerful tools for addressing these challenges: the regularity structure theory developed by Hairer \cite{Hairer14}, and the paracontrolled distributions approach introduced by Gubinelli, Imkeller, and Perkowski \cite{GIP15}. Both frameworks leverage the intrinsic structure of solutions to provide rigorous interpretations for terms that are not classically well-defined. In \cite{HZZZ24}, the authors introduced the paracontrolled distribution method and performed renormalization by probabilistic calculation to establish the global well-posedness of \eqref{Koleq} for $\a\in(-\ff{2}{3},-\ff{1}{2})$. However, to the best of our knowledge, the case where $\a\in(-1,-\ff{2}{3}]$ has not yet been explored in the existing literature. In this paper, we aim to investigate this previously unexamined case. Furthermore, we introduce weighted spaces to derive uniform bounds for the solutions of the linear equation with a polynomial growth, while also handling cases where $b$ includes certain Gaussian random fields (see Example \ref{eg5.3}). This example has been attracting growing research interest due to its close connection with quantum field theory, in particular the stochastic quantization raised by Parisi and Wu \cite{PW81}.

On the other hand, recently, in \cite{HZ23}, Hao and Zhang investigated the following first-order stochastic differential equation (SDE) with distributional drifts:
\begin{align}\label{nondege}
\dif X_t = b(t,X_t)\dif t + \sqrt{2}\dif W_t, 
\end{align}
and established the well-posedness of SDE \eqref{nondege} for both subcritical and supercritical cases. In the subcritical case, they applied Zvonkin's transformation, whereas in the supercritical case, they analyzed the generalized martingale problem by using an energy method. Motivated by this work, we aim to extend the results to the kinetic case. Notably, we must emphasize that the equation \eqref{3} is degenerate, rendering Zvonkin's transformation potentially inapplicable, which brings us great difficulties. Roughly speaking, when $\a\leq0$, $u\in\bC_\bba^{2+\a}\subset\bC_x^{\ff{2+\a}{3}}$ with $\ff{2+\a}{3}<1$, which implies that $\nn_x u$ may not exist. To tackle this problem, we establish a  generalized version of  It\^{o}'s formula (see Lemma \ref{geito}) and consider the associated martingale problem. Thus far, we have only been able to handle the subcritical case. The challenge in the supercritical case lies in dealing with the kinetic Kolmogorov equations with supercritical drifts, which remains a target for future research.

\subsection{Structure of the paper}

This paper is organized as follows:

In Section \ref{sec21} we introduce the notion of anisotropic H\"{o}lder space. To make $b \cdot \nabla_v u$ well-defined in PDE \eqref{Koleq} when $\alpha_b < -\ff{1}{2}$, paraproduct calculus is introduced in Section \ref{sec22}.

Section \ref{sec3} focuses on the solution of equation \eqref{Koleq} with a weighted distributional drift $b$. Under the condition {\bf (H$_\text{w}^\text{sub}$)} below, we obtain the unique solution to equation \eqref{Koleq}. The key technique is the localization, which is used in \cite{HZZZ24}.

Section \ref{sec4} is dedicated to the study of SDEs with subcritical distributional drifts, where we prove Theorem \ref{sec1:main2}. In Subsection \ref{sec41}, we introduce the notion of Krylov's estimate and define a generalized form of Young’s integral. We emphasize that the Krylov's estimate guarantees that the process $A_t^b$ is an energy-zero process, enabling us to apply the generalized It\^{o} formula to prove the uniqueness of the solution. Then we establish a crucial uniform Krylov's estimate for the approximated solutions $Z^n$ under {\bf (H$^\text{sub}$)} and {\bf (H$_\text{w}^\text{sub}$)}, respectively. In Subsection \ref{sec42}, we follow a standard approach to show the existence of the solutions to SDE \eqref{3} by utilizing Prokhorov’s theorem and Skorokhod’s representation theorem. In Subsection \ref{sec43}. we apply the generalized It\^o's formula and study the associated martingale problem to prove the uniqueness.

Section \ref{sec5} presents an application of Theorem \ref{sec1:main2} by studying the drift of Gaussian random fields. 
We give the conditions for the Gaussian random fields $b$ under which SDE \eqref{3} is well-posed.

We conclude the introduction with the following notations. Throughout this paper, we use $C$ with or without subscripts to denote constants, whose values may vary from line to line. Also, we use $:=$ to indicate a definition. We write $A\lesssim_C B$ and $A\asymp_C B$ or simply without subscripts to mean that for some constant $C\geq1$, $A\leq CB$ and $C^{-1}B\leq A\leq CB$. Moreover,
for the convenience of the readers, we list some frequently used notations:

\begin{itemize}
\item $\mN$: set of positive integers; $\mN_0:=\{0\}\cup\mN$.

\item $\mR^n$: Euclidean space of dimension $n\in\mN$.

\item $\bC^s_\bba(\rho_\k)$: weighted anisotropic H\"{o}lder space.

\item $P_t$: semigroup of the operator $\Delta_v+v\cdot\nn_x$, which is given by 
$$P_t f=\G_tp_t*\G_tf=\G_t(p_t*f),$$
where $\Gamma_t f(z):=f(\Gamma_t z):=f(x-tv,v)$.

\item $f\prec g, f\circ g, f\succ g$: paraproduct. $f\preceq g:=f\prec g+f\circ g$ (see \eqref{bonyeq} below).

\end{itemize}

\section{Preliminaries}\label{sec2}

In this section we introduce the basic notations and propositions of weighted anisotropic H\"{o}lder spaces (cf. \cite{HZZZ24} and \cite{Triebel06}).
Let $\sS(\mR^{2d})$ be the Schwarz space of all rapidly decreasing functions on $\mR^{2d}$, and let $\sS'(\mR^{2d})$ denote the dual space of $\sS(\mR^{2d})$, known as Schwartz generalized function (or tempered distribution) space. For any $f\in\sS(\mR^{2d})$, we define the Fourier transform $\hat f$ and inverse Fourier transform $\check f$ respectively by
$$\hat f(\xi):=\ff{1}{(2\pi)^{d}}\int_{\mR^{2d}}\e^{-i\xi\cdot z}f(z)\dif z,\ \xi\in\mR^{2d},$$
$$\check f(z):=\ff{1}{(2\pi)^{d}}\int_{\mR^{2d}}\e^{i\xi\cdot z}f(\xi)\dif\xi, \ z\in\mR^{2d}.$$
For $q\in[1,\infty]$ and a normed space $\mB$, we will use the notation
$$\mL_T^q\mB:=L^q([0,T];\mB),\ \mL_T^q:=\mL_T^qL^q.$$
Additionally, we denote by $C_b=C_b(\mR^{2d})$ (resp. $C_c=C_c(\mR^d)$) the Banach space of all bounded continuous functions on $\mR^{2d}$ (resp. the space of all continuous functions on $\mR^{2d}$ with compact support), and by $C_b^\infty=C_b^\infty(\mR^{2d})$ (resp. $C_c^\infty=C_c^\infty(\mR^{2d})$) the space of all smooth functions with bounded derivatives of all orders (resp. the space of all smooth functions with compact support). Furthermore, we define the weighted space $C_b(\rho_\k):=\{f:f\rho_\k\in C_b\}$.

\subsection{Anisotropic weighted H\"{o}lder spaces}\label{sec21}

In this subsection we recall the definition of weighted anisotropic H\"{o}lder space.
For $r>0$ and $z\in\mR^{2d}$, the ball centered at $z$ and with radius $r$ in terms of the anisotropic distance 
$|\cdot|_\bba$ (see \eqref{distance}) is defined as follows:
$$B_r^\bba(z):=\{z'\in\mR^{2d}:|z-z'|_\bba\leq r\},\ B_r^\bba:=B_r^\bba(0).$$
Let $\chi_0^\bba$ be a symmetric $C^\infty$-function on $\mR^{2d}$ such that
$$\chi_0^\bba(\xi)=1 \text{\ for\ } \xi\in B_1^\bba \text{\ and\ }\chi_0^\bba(\xi)=0 \text{\ for\ } \xi\notin B_{4/3}^\bba.$$
For $j\in\mN$, define 
\begin{align}
\phi_j^\bba(\xi):=\left\{
\begin{aligned}
&\chi_0^\bba(2^{-j\bba}\xi)-\chi_0^\bba(2^{-(j-1)\bba}\xi),\ \ &j\geq1,\\
&\chi_0^\bba(\xi),&j=0,
\end{aligned}
\right.\nonumber
\end{align}
where for $s\in\mR$ and $\xi=(\xi_1,\xi_2)$, $2^{s\bba}\xi:=(2^{3s}\xi_1,2^s\xi_2).$ This definition ensures that 
\begin{align}\label{2.1}
\sum_{j\geq0}\phi_j^\bba(\xi)=1,\ \forall\xi\in\mR^{2d},
\end{align}
and
\begin{align}\label{supp}
\supp(\phi_j^\bba)\subset\{\xi:2^{j-1}\leq|\xi|_\bba\leq2^{j+1}\},\ j\geq1,\ \supp(\phi_0^\bba)\subset B_{4/3}^\bba.
\end{align}
For given $j\geq0$, we define the dyadic anisotropic block operator $\cR_j^\bba$ on $\sS'$ as
\begin{align}
\cR_j^\bba f(z):=(\phi_j^\bba\hat f)\check\ (z)=\check\phi_j^\bba*f(z),
\end{align}
where the convolution is taken in the distributional sense. By scaling, we have
\begin{align}\label{2.3}
\check\phi_j^\bba(z)=2^{4d(j-1)}\check\phi_1^\bba(2^{(j-1)\bba}z),\ j\geq1.
\end{align}
For $j\in\mN$, it follows that
\begin{align}\label{wtR}
\cR_j^\bba=\cR_j^\bba\wt\cR_j^\bba,\ \text{where}\ \wt\cR_j^\bba:=\sum_{|i-j|\leq2}\cR_i^\bba,
\end{align}
with the convention that $\cR_j^\bba:=0$ for $j<0$.

The following definition about the admissible weights is from \cite{Triebel06}.
\bd
A $C^\infty$-smooth function $\rho:\mR^{2d}\to (0,\infty)$ is called an admissible weight if for each $j\in\mN$, there is a constant $C_j>0$ such that
\begin{align*}
    |\nabla^j \rho(z)|\le C_j \rho(z),\ \forall z\in\mR^{2d},
\end{align*}
and for some $C,\kappa>0$,
\begin{align*}
    \rho(z)\le C\rho(z')(1+|z-z'|_\bba^\kappa),\ \forall z,z'\in\mR^{2d}.
\end{align*}
For any admissible weight $\rho$ and $p\in[1,\infty]$, we define $\|f\|_{L^p(\rho)}:=\|\rho f\|_{L^p}$.
\ed

\begin{remark}
It can be readily verified that for any $\k\in\mR$, the weight function $ \rho_\k$ defined in \eqref{rhok} is an admissible weight.
\end{remark}

The following Bernstein's inequality is standard and can be found in \cite[Lemma 2.4]{HZZZ24}.

\bl\label{bernstein}
Let $\rho$ be an admissible weight. For any $\bbk=(k_1,k_2)\in\mN^2$ and $1\leq p\leq q\leq\infty$, there exists a constant $C=C(\rho,\bbk,p,q,d)>0$ such that for all $j\geq0$, 
\begin{align}\label{bernstein2}
\|\nn_x^{k_1}\nn_v^{k_2}\cR_j^\bba f\|_{L^q(\rho)}\lesssim_C 2^{j(\bba\cdot\bbk+4d(\ff{1}{p}-\ff{1}{q}))}\|\cR_j^\bba f\|_{L^p(\rho)}.
\end{align}
\el

Now we define the weighted anisotropic H\"{o}lder spaces as follows (see \cite{Triebel06, HZZZ24}).

\bd\label{Holder}
Let $\rho$ be an admissible weight. For any $s\in\mR\backslash\mN_0$, the anisotropic H\"{o}lder space is defined by
\begin{align}
\bC_\bba^s(\rho):=\Big\{f\in\sS':\|f\|_{\bC_\bba^s(\rho)}:=\sup_{j\geq0}2^{js}\|\cR_j^\bba f\|_{L^\infty(\rho)}<\infty\Big\}.\no
\end{align}
For $s\in\mN_0$, the anisotropic H\"{o}lder space is defined by
\begin{align}
\bC_\bba^s(\rho):=\Big\{f\in\sS':\|f\|_{\bC_\bba^s(\rho)}:=\sup_{j\geq0}2^{js}\|\cR_j^\bba f\|_{L^\infty(\rho)}+\|\nn_v^s f\|_{L^\infty(\rho)}<\infty\Big\}\no
\end{align}
\ed

We have the following interpolation inequality.

\bl
For any $s_0,s_1\in\mR$ with $s_0\neq s_1$ and $\theta\in(0,1)$, there exists a constant $C=C(d,s_0,s_1,\theta)>0$ such that for any admissible weight $\rho$ and $f\in\bC_\bba^{s_0\vee s_1}(\rho)$,
\begin{align}\label{interpolation}
\|f\|_{\bC_\bba^{s_\theta}(\rho)}\lesssim_C \|f\|_{\bC_\bba^{s_0}(\rho)}^{\theta}\|f\|_{\bC_\bba^{s_1}(\rho)}^{1-\theta},
\end{align}
where $s_\theta:=\theta s_0+(1-\theta) s_1$.
\el

\bpf
By \cite[Theorem 2.80]{BCD11}, we have for $s_\theta\in\mR\backslash\mN_0$,
\begin{align}
\|f\|_{\bC_\bba^{s_\theta}(\rho)}\overset{\text{\cite[(2.22)]{HZZZ24}}}{=}\|\rho f\|_{\bC_\bba^{s_\theta}}\lesssim \|\rho f\|_{\bC_\bba^{s_0}}^{\theta}\|\rho f\|_{\bC_\bba^{s_1}}^{1-\theta}=\|f\|_{\bC_\bba^{s_0}(\rho)}^{\theta}\|f\|_{\bC_\bba^{s_1}(\rho)}^{1-\theta}.\no
\end{align}
For $s_\theta\in\mN_0$, similarly
\begin{align}
&\|f\|_{\bC_\bba^{s_\theta}(\rho)}\lesssim\sup_{j\geq0}2^{js_\theta}\|\cR_j^\bba f\|_{L^\infty(\rho)}+\sum_{j\geq0}\|\nn_v^{s_\theta} \cR_j^\bba f\|_{L^\infty(\rho)}\overset{\eqref{bernstein2}}{\lesssim}\sum_{j\geq0}2^{js_\t}\|\cR_j^\bba f\|_{L^\infty(\rho)}\no\\
&\qquad\quad\overset{\text{\cite[(2.22)]{HZZZ24}}}{=}\sum_{j\geq0}2^{js_\t}\|\cR_j^\bba (\rho f)\|_{L^\infty}\lesssim \|\rho f\|_{\bC_\bba^{s_0}}^{\theta}\|\rho f\|_{\bC_\bba^{s_1}}^{1-\theta}=\|f\|_{\bC_\bba^{s_0}(\rho)}^{\theta}\|f\|_{\bC_\bba^{s_1}(\rho)}^{1-\theta}.\no
\end{align}
The proof is complete.
\epf

\subsection{Paraproduct calculus}\label{sec22}
To ensure that the term $b\cdot\nn_v u$ makes sense when $b$ is distribution-valued and $u\in\sS(\mR^{2d})$, we recall the well-known Bony's decomposition, which allows us to control the product of functions in anisotropic H\"{o}lder spaces under weaker regularity assumptions. For any $k\geq0$, define the cut-off low frequency operator $S_k^\bba$ as
$$S_k^\bba f:=\sum_{j=0}^{k-1}\cR_j^\bba f\overset{k\rightarrow\infty}{\rightarrow}f,$$
where we take $S_{-1}^\bba=0$ by convention. For $f,g\in\sS'(\mR^{2d})$, we define the following paraproducts
$$f\prec g:=\sum_{k\geq0}S_{k-1}^\bba f\cR_k^\bba g,\ f\circ g:=\sum_{k\geq0}\cR_k^\bba f\wt\cR_k^\bba g,$$
where $\wt\cR_j^\bba$ is defined in \eqref{wtR}. The Bony's decomposition of $fg$ is formally written as (cf. \cite{BCD11})
\begin{align}\label{bonyeq}
fg=f\prec g+f\circ g+g\prec f=:f\preceq g+f\succ g.
\end{align} 
The essential property of the Bony's decomposition is that
\begin{align}
\cR_k^\bba(S^\bba_{j-1}f\cR_j^\bba g)=0\ \text{for}\ |k-j|>2.
\end{align}

We now state the key paraproduct estimates (cf. \cite[Lemma 2.11]{HZZZ24}).

\bl\label{bony}
Let $\rho_1,\rho_2$ be admissible weights and $s_1, s_2\in\mR$ with $s_1+s_2\notin\mN_0$.
\begin{enumerate}[]
\item {\rm\text{(i)}} If $s_1<0$, then there exists a constant $C=C(d,s_1,s_2)>0$ such that
\begin{align}\label{paraproductcontrol1}
\| f\prec g\|_{\bC_\bba^{s_1+s_2}(\rho_1\rho_2)}\lesssim_C\|f\|_{\bC_\bba^{s_1}(\rho_1)}\|g\|_{\bC_\bba^{s_2}(\rho_2)}.
\end{align}

For $s_1=0$, we have 
\begin{align}\label{paraproductcontrol2}
\| f\prec g\|_{\bC_\bba^{s_2}(\rho_1\rho_2)}\lesssim_C\|f\|_{L^\infty(\rho_1)}\|g\|_{\bC_\bba^{s_2}(\rho_2)}.
\end{align}

\item {\rm\text{(ii)}} If $s_1+s_2>0$, then there exists a constant $C=C(d,s_1,s_2)>0$ such that
\begin{align}\label{paraproductcontrol3}
\|f\circ g\|_{\bC_\bba^{s_1+s_2}(\rho_1\rho_2)}\lesssim_C\|f\|_{\bC_\bba^{s_1}(\rho_1)}\|g\|_{\bC_\bba^{s_2}(\rho_2)}.
\end{align}
\end{enumerate}
\el

\br
By the above lemma, if $s_1+s_2>0$ with $s_1+s_2\notin\mN_0$, then there exists a constant $C=C(d,s_1,s_2)>0$ such that
\begin{align}\label{paraproductcontrol4}
\|fg\|_{\bC_\bba^{s_1\wedge s_2}(\rho_1\rho_2)}\lesssim_C \|f\|_{\bC_\bba^{s_1}(\rho_1)}\|g\|_{\bC_\bba^{s_2}(\rho_2)}.
\end{align}
\er

Using Bony's decomposition we can write
\begin{align}\label{decomposition}
b\cdot\nabla_v u:={\div_v(b\preceq u)+b\succ\nabla_v u}-\div_v b\preceq u=:b\odot\nabla_v u-\div_v b\preceq u.
\end{align}
In particular, if $\div_v b=0$, then 
$$b\cdot\nabla_v u=b\odot\nabla_v u.$$

The following regularity estimates are crucial in establishing Theorem \ref{pdees}.

\bl\label{driftbony}
Let $\rho_1,\rho_2$ be admissible weights. For any $s\in(-1,0)$, there exists a constant $C=C(d,s)>0$ such that 
\begin{align}\label{bonycontrol1}
\|b\odot\nabla_v u\|_{\bC_\bba^s(\rho_1\rho_2)}\lesssim_C\|b\|_{\bC_\bba^s(\rho_1)}\|u\|_{\bC_\bba^1(\rho_2)}
\end{align}
and
\begin{align}\label{bonycontrol2}
\|\div_v b\preceq u\|_{\bC_\bba^{1+s}(\rho_1\rho_2)}\lesssim_C\|\div_v b\|_{\bC_\bba^s(\rho_1)}\|u\|_{\bC_\bba^1(\rho_2)}.
\end{align}
\el

\bpf
For $s\in(-1,0)$, by Lemma \ref{bernstein} and Lemma \ref{bony} we have
\begin{align}
\|\div_v(b\preceq u)\|_{\bC_\bba^s(\rho_1\rho_2)}\lesssim\|b\preceq u\|_{\bC_\bba^{s+1}(\rho_1\rho_2)}\lesssim\|b\|_{\bC_\bba^s(\rho_1)}\|u\|_{\bC_\bba^1(\rho_2)},\no
\end{align}
and 
\begin{align}
\|b\succ\nabla_v u\|_{\bC_\bba^s(\rho_1\rho_2)}\lesssim\|b\|_{\bC_\bba^s(\rho_1)}\|\nn_v u\|_{L^\infty(\rho_2)}\lesssim\|b\|_{\bC_\bba^s(\rho_1)}\|u\|_{\bC_\bba^{1}(\rho_2)}.\no
\end{align}
Estimate \eqref{bonycontrol1} then follows directly from the above bounds,
while \eqref{bonycontrol2} can be obtained analogously via a similar calculation.
\epf

\section{PDEs with weighted distributional drifts}\label{sec3}

In this section we study the solvability of the following PDE with $\l\geq0$:
\begin{align}\label{weighted distributional PDE}
\p_t u=\Delta_v u-v\cdot\nn_x u-\l u+b\cdot\nn_v u+f,\ u_0\equiv0.
\end{align}
Below we always assume that the following condition on $b$ holds:
\begin{enumerate}[{\bf (H$_\text{w}^\text{sub}$)}]
\item Let $\a_b\in(-1,0)$, $q_b\in(\ff{2}{1+\a_b},\infty]$ and $\k\in[0,\ff{1+\a_b-2/q_b}{3+\a_b-2/q_b}]$. Suppose that
\begin{align}\label{c_b}
c_b:=\|b\|_{\mL_T^{q_b}\bC^{\a_b}_\bba(\rho_\k)}+\|\div_v b\|_{\mL_T^{q_b}\bC^{\a_b}_{\bba}}<\infty.
\end{align}
\end{enumerate}
For convenience of notations, we introduce the following parameter set:
$$\Xi:=(T,d,\a_b,q_b,\k,c_b).$$

Let $p_t(z)=p_t(x,v)$ be the distributional density of $Z_t:=(\ss{2}\int_0^t B_s\dif s, \ss{2}B_t)$, given by
$$p_t(x,v)=\Big(\ff{2\pi t^4}{3}\Big)^{-\ff{d}{2}}\exp\Big(-\ff{3|x|^2+|3x-2tv|^2}{4t^3}\Big).$$
It is easy to see that for any $\l>0$,  
\begin{align}\label{pscaling}
p_{\l t}(x,v)=\l^{-2d}p_t(\l^{-3/2}x,\l^{-1/2}v).
\end{align} 
For $t>0$, let $P_t$ be the kinetic semigroup defined by
\begin{align}\label{semigroup}
P_tf(z):=\mE f(\Gamma_t z+Z_t)=\Gamma_t(p_t*f)(z)=(\Gamma_t p_t*\Gamma_t f)(z),
\end{align}
where 
$$\Gamma_t f(z):=f(\Gamma_t z):=f(x-tv,v).$$
For $\l\geq0$ and a time-dependent distribution $f_t(\cdot):\mR_+\to\sS'(\mR^{2d})$, we define
\begin{align}\label{duhamel}
u_t(z):=\sI_t^\l(f)(z):=\int_0^t \e^{-\l(t-s)}P_{t-s}f_s(z)\dif s,\ t>0.
\end{align}
Then it is easy to see that in the distributional sense,
\begin{align}\label{2.19}
\partial_t u=(\Delta_v-v\cdot\nn_x-\l)u+f.
\end{align}
Next we give the definition of the solution to \eqref{weighted distributional PDE}.

\bd\label{mildsol}
Let $T>0$ and $\k\in\mR$. We call  $u\in\mL_T^\infty\bC_\bba^1(\rho_\k)$ a mild solution of PDE \eqref{weighted distributional PDE}, if $u$ solves the following integral equation:
$$u_t=\sI_t^\l(b\odot\nabla_v u-\div_v b\preceq u+f).$$
\ed

\br\label{remark3.3}
Under {\bf (H$_\text{w}^\text{sub}$)}, by Lemma \ref{driftbony} one sees that all terms on the right side of the equality are well-defined. By \eqref{2.19}, the solution $u$ constructed above satisfies that for any test function $\varphi\in C_c^\infty(\mR^{2d})$, 
$$\p_t\<u,\varphi\>=\<u,\Delta_v\varphi\>+\<u,v\cdot\nn_x\varphi\>-\l\<u,\varphi\>+\<b\odot\nabla_v u-\div_v b\preceq u,\varphi\>+\<f,\varphi\>,$$
where $\<u,\varphi\>:=\int_{\mR^{2d}}u(z)\varphi(z)\dif z$. This formulation indicates that $u$ is also a weak solution of PDE \eqref{weighted distributional PDE}.
\er

The following a priori estimate with unweighted drifts is standard.
\bl\label{l3.1}
Assume that $\a_b\in(-1,0)$, $q_b\in(\ff{2}{1+\a_b},\infty]$ and 
\begin{align}\label{ell^b}
\ell_b:=\|b\|_{\mL_T^{q_b}\bC^{\a_b}_{\bba}}+\|\div_v b\|_{\mL_T^{q_b}\bC^{\a_b}_{\bba}}<\infty.
\end{align}
Then for any $q\in(\ff{2}{1+\a_b},q_b]$, there is a constant $C=C(T,d,\a_b,q_b,\ell_b)>0$ such that for any solution $u$ to the PDE \eqref{weighted distributional PDE},
\begin{align}\label{6.6}
\|u\|_{\mL_T^\infty\bC_\bba^{2+\a_b-2/q}}\lesssim_C\|f\|_{\mL_T^q\bC_\bba^{\a_b}}.
\end{align}
\el

\bpf
Let $q_0:=\ff{q_bq}{q_b-q}$. By Lemma \ref{driftbony} and H\"older's inequality, we have
\begin{align}\label{bonybnnu}
\|b\cdot\nn_v u\|_{\mL_T^q\bC_\bba^{\a_b}}\overset{\eqref{decomposition}}{\leq}\|b\odot\nn_v u\|_{\mL_T^q\bC_\bba^{\a_b}}+\|\div_v b\preceq u\|_{\mL_T^q\bC_\bba^0}\lesssim\ell_b\|u\|_{\mL_T^{q_0}\bC_\bba^1}.
\end{align}
By \cite[Theorem 2.17]{HRZ23}, we have the following Schauder's estimate:
\begin{align}
\|\sI_.^\l(f)\|_{\mL_T^{\infty}\bC_\bba^{2+\a_b-2/q}}\lesssim_C\|f\|_{\mL_T^q\bC_\bba^{\a_b}}.\no
\end{align}
Thus,
\begin{align}\label{pf:le3.3:1}
\|u\|_{\mL_T^\infty\bC_\bba^{2+\a_b-2/q}}&\lesssim\|b\cdot\nn_v u\|_{\mL_T^q\bC_\bba^{\a_b}}+\|f\|_{\mL_T^q\bC_\bba^{\a_b}}\lesssim\ell_b\|u\|_{\mL_T^{q_0}\bC_\bba^1}+\|f\|_{\mL_T^q\bC_\bba^{\a_b}},
\end{align}
which implies that for any $t\in[0,T]$ and $q\in(\ff{2}{1+\a_b},q_b)$, 
\begin{align}
\|u(t)\|^{q_0}_{\bC_\bba^{2+\a_b-2/q}}\lesssim\ell_b\int_0^t \|u(s)\|^{q_0}_{\bC_\bba^{2+\a_b-2/q}}\dif s+\|f\|^{q_0}_{\mL_T^q\bC_\bba^{\a_b}}.\no
\end{align}
Hence, by Gronwall's inequality,
\begin{align}
\|u\|_{\mL_T^\infty\bC_\bba^{2+\a_b-2/q}}\lesssim\|f\|_{\mL_T^q\bC_\bba^{\a_b}}.\no
\end{align}
Note that for $q\in(\ff{2}{1+\a_b},q_b)$,
\begin{align}\label{pf:le3.3:2}
\|u\|_{\mL_T^{q_0}\bC_\bba^1}\lesssim\|u\|_{\mL_T^\infty\bC_\bba^{2+\a_b-2/q}}\lesssim\|f\|_{\mL_T^q\bC_\bba^{\a_b}}.
\end{align}
Substituting \eqref{pf:le3.3:2} into \eqref{pf:le3.3:1}, we conclude that \eqref{6.6} still holds for $q=q_b$.
\epf

To establish the main result of this section, we need the following theorem, which provides a sharper estimate than Lemma \ref{l3.1}. More concretely, the right side of \eqref{6.7} depends polynomially on $\ell_b$, rather than exponentially as in Lemma \ref{l3.1}. 

\bt\label{theorem3.2}
Assume that $\a_b\in(-1,0)$, $q_b\in(\ff{2}{1+\a_b},\infty]$ and
\begin{align}
\ell_b:=\|b\|_{\mL_T^{q_b}\bC^{\a_b}_{\bba}}+\|\div_v b\|_{\mL_T^{q_b}\bC^{\a_b}_{\bba}}<\infty.\no
\end{align}
Let $\a\in(\ff{2}{q_b}-1,\a_b]$ and $q\in(\ff{2}{1+\a},q_b]$. Then for any $f\in\mL_T^q\bC_\bba^\a$ and $\l\geq0$, there is a unique solution $u$ to the PDE \eqref{weighted distributional PDE} such that for any $\t\in[0,2-\ff{2}{q})$,
\begin{align}\label{6.7}
(1+\l)^{\ff{\t}{2}}\|u\|_{\mL_T^\infty\bC_\bba^{1+\b(\a,q)-\t}}\lesssim_C (1+\ell_b)^{\ff{2}{\b(\a,q)}}\|f\|_{\mL_T^q\bC_\bba^\a},
\end{align}
and
\begin{align}\label{6.8}
\|u\|_{\mL_T^\infty}\lesssim_C (1+\ell_b)^{\ff{1-\b(\a,q)}{\b(\a,q)}}\|f\|_{\mL_T^q\bC_\bba^\a},
\end{align}
where $\b(\a,q):=1+\a-\ff{2}{q}>0$ and the constant $C:=C(\Xi,\delta,\t,q)>0$.
\et

To prove this theorem, we need the following maximum estimate for solutions in the sense of Definition \ref{mildsol}.

\bl
Assume that $\a_b\in(-1,0)$, $q_b\in(\ff{2}{1+\a_b},\infty]$ and
\begin{align}
\ell_b:=\|b\|_{\mL_T^{q_b}\bC^{\a_b}_{\bba}}+\|\div_v b\|_{\mL_T^{q_b}\bC^{\a_b}_{\bba}}<\infty.\no
\end{align}
Then for any $f\in\mL_T^\infty$, there is a unique solution $u\in\mL_T^\infty\bC_{\bba}^1$ to the PDE \eqref{weighted distributional PDE} such that
\begin{align}\label{6.4}
\|u\|_{\mL_T^\infty}\leq T\|f\|_{\mL_T^\infty}.
\end{align}
\el

\bpf
Let $b_n(t,z):=b(t,\cdot)*\G_n(z)$. It is well-known that there exists a unique solution $u_n\in\mL_T^\infty C_b^2$ to the PDE \eqref{weighted distributional PDE} with only replacing $b$ by $b_n$ (cf. \cite[Theorem 1.3]{CHM21}). By \cite[Theorem 6.2]{HWZ20}, we have
\begin{align}\label{5.4}
\|u_n\|_{\mL_T^\infty}\leq T\|f\|_{\mL_T^\infty}.
\end{align}
For $\a_b\in(-1,0)$ and $t\in[0,T]$, by definition we have
\begin{align}
\|b_n(t)\|_{\bC_\bba^{\a_b}}&=\sup_{j\geq0}2^{j\a_b}\|\cR_j^\bba b_n(t)\|_{L^\infty}=\sup_{j\geq0}2^{j\a_b}\|\check\phi_j^\bba*b(t)*\Gamma_n\|_{L^\infty}\no\\
&\leq\sup_{j\geq0}2^{j\a_b}\|\check\phi_j^\bba*b(t)\|_{L^\infty}\|\Gamma_n\|_{L^1}=\sup_{j\geq0}2^{j\a_b}\|\cR_j^\bba b(t)\|_{L^\infty}=\|b(t)\|_{\bC_\bba^{\a_b}}.\no
\end{align}
By a similar calculation for $\div_v b$, we obtain
\begin{align}\label{molifycompare}
\ell_{b_n}=\|b_n\|_{\mL_T^{q_b}\bC_\bba^{\a_b}}+\|\div_v b_n\|_{\mL_T^{q_b}\bC_\bba^{\a_b}}\leq\|b\|_{\mL_T^{q_b}\bC_\bba^{\a_b}}+\|\div_v b\|_{\mL_T^{q_b}\bC_\bba^{\a_b}}=\ell_{b}.
\end{align}
Then by Lemma \ref{l3.1} we have the following uniform estimates: for any $q\in(\ff{2}{1+\a_b},q_b]$,
\begin{align}\label{3.14}
\sup_n\|u_n\|_{\mL_T^\infty\bC_\bba^{2+\a_b-2/q}}\lesssim\|f\|_{\mL_T^q\bC_\bba^{\a_b}}\lesssim\|f\|_{\mL_T^\infty}.
\end{align}
By a classical compactness argument based on the Aubin-Lions lemma (see \cite[Lemma 6.9]{Zhang21}), there exists a  function $u\in\mL_T^\infty\bC_\bba^1$ and a subsequence $\{n_k\}$ such that 
\begin{align}
\lim_{k\to\infty}\|u_{n_k}-u\|_{\mL_T^\infty\bC_\bba^1}.\no
\end{align}
Since the PDE \eqref{weighted distributional PDE} is linear, a standard argument (cf. \cite[Theorem 3.6]{HZ23}) yields that $u$ is the unique weak solution to \eqref{weighted distributional PDE}, and \eqref{6.4} follows by \eqref{5.4}. This completes the proof.
\epf

Now we can give

\bpf[Proof of Theorem \ref{theorem3.2}]
We only prove the a priori estimates \eqref{6.7} and \eqref{6.8}. By \cite[Theorem 2.17]{HRZ23}, for $\b(\a,q)=1+\a-\ff{2}{q}$ and $\t\in[0,2-\ff{2}{q})$ we have
\begin{align}\label{pf:thm3.2:1}
(1+\l)^{\ff{\t}{2}}\|u\|_{\mL_T^\infty\bC_\bba^{1+\b(\a,q)-\t}}\lesssim\|b\cdot\nn_v u\|_{\mL_T^q\bC_\bba^{\a_b}}+\|f\|_{\mL_T^q\bC^\a_\bba}\overset{\eqref{bonybnnu}}{\lesssim}\ell_b\|u\|_{\mL_T^{q_0}\bC_\bba^1}+\|f\|_{\mL_T^q\bC_\bba^\a},
\end{align}
where $q_0:=\ff{q_bq}{q_b-q}$. In particular, for $\t=\b(\a,q)$, there exist constants $C_i=C_i(\a_b,\a,d,T,q)>0$, $i=1,2$, such that for all $\l\geq0$,
$$(1+\l)^{\ff{\b(\a,q)}{2}}\|u\|_{\mL_T^\infty\bC_\bba^1}\leq C_1\ell_b\|u\|_{\mL_T^{q_0}\bC_\bba^1}+C_2\|f\|_{\mL_T^q\bC_\bba^\a}.$$
Define $\l_b:=(2C_1\ell_b)^{\ff{2}{\b(\a,q)}}$. For all $\l\geq\l_b$ we have
\begin{align}
2C_1\ell_b\|u\|_{\mL_T^\infty\bC_\bba^1}\leq(1+\l)^{\ff{\b(\a,q)}{2}}\|u\|_{\mL_T^\infty\bC_\bba^1}\leq C_1\ell_b\|u\|_{\mL_T^{q_0}\bC_\bba^1}+C_2\|f\|_{\mL_T^q\bC_\bba^\a},\no
\end{align}
which implies that
$$\ell_b\|u\|_{\mL_T^\infty\bC_\bba^1}\lesssim\|f\|_{\mL_T^q\bC_\bba^\a}.$$
Substituting this into \eqref{pf:thm3.2:1}, we get for any $\t\in[0,2-\ff{2}{q})$,
\begin{align}\label{6.9}
(1+\l)^{\ff{\t}{2}}\|u\|_{\mL_T^\infty\bC_\bba^{1+\b(\a,q)-\t}}\lesssim\|f\|_{\mL_T^q\bC_\bba^\a}.
\end{align}
To drop the condition $\l\geq\l_b$, we make the following decomposition $u=u_1+u_2$:
\begin{align}
&\p_t u_1=\Delta_v u_1-v\cdot\nn_x u_1-(\l+\l_b) u_1+b\cdot\nn_v u_1+f,\ u_1(0)=0,\no\\
&\p_t u_2=\Delta_v u_2-v\cdot\nn_x u_2-\l u_2+b\cdot\nn_v u_2+\l_b u_1,\ u_2(0)=0.\no
\end{align}
Taking $\t=1+\b(\a,q)$ in \eqref{6.9} for $u_1$, we have
\begin{align}
\|u\|_{\mL_T^\infty}&\leq\|u_1\|_{\mL_T^\infty}+\|u_2\|_{\mL_T^\infty}\overset{\eqref{6.4}}{\leq}(1+\l_b T)\|u_1\|_{\mL_T^\infty}\no\\
&\lesssim(1+\l_b)(1+\l_b)^{-\ff{1+\b(\a,q)}{2}}\|f\|_{\mL_T^q\bC_\bba^\a}\lesssim(1+\ell_b)^{\ff{1-\b(\a,q)}{\b(\a,q)}}\|f\|_{\mL_T^q\bC_\bba^\a},
\end{align}
which is \eqref{6.8}. To prove \eqref{6.7}, by \eqref{interpolation}  with unweighted interpolation and Young's inequality, we have for any $\t\in[0,2-\ff{2}{q})$ and $\l\geq0$, there exist $C_i=C_i(\a_b,\a,d,T,q,\t)>0$, $i=3,4$ and $\varepsilon>0$ such that
\begin{align}\label{6.11}
&(1+\l)^{\ff{\t}{2}}\|u\|_{\mL_T^\infty\bC_\bba^{1+\b(\a,q)-\t}} \overset{\eqref{pf:thm3.2:1}}{\leq} C_3\ell_b\|u\|_{\mL_T^{q_0}\bC_\bba^1}+C_4\|f\|_{\mL_T^q\bC_\bba^\a}\no\\
&\qquad\leq C_3\ell_b\|u\|_{\mL_T^{q_0/(1+\b(\a,q))} \bC_\bba^{1+\b(\a,q)}}^{\ff{1}{1+\b(\a,q)}}\|u\|_{\mL_T^\infty}^{\ff{\b(\a,q)}{1+\b(\a,q)}}+C_4\|f\|_{\mL_T^q\bC_\bba^\a}\no\\
&\qquad\leq C_\varepsilon(C_3\ell_b)^{\ff{1+\b(\a,q)}{\b(\a,q)}}\|u\|_{\mL_T^\infty}+\varepsilon\|u\|_{\mL_T^{q_0/(1+\b(\a,q))} \bC_\bba^{1+\b(\a,q)}}+C_4\|f\|_{\mL_T^q\bC_\bba^\a}.
\end{align}
Letting $\t=0$ and $\varepsilon=\ff{1}{2}$, by \eqref{6.8} we have 
$$\|u\|_{\mL_T^\infty\bC_\bba^{1+\b(\a,q)}}\lesssim(1+\ell_b)^{\ff{2}{\b(\a,q)}}\|f\|_{\mL_T^q\bC_\bba^\a}.$$
Substituting this into \eqref{6.11}, we get \eqref{6.7}.
\epf

Now we prove the main result of this section.

\bt\label{pdees}
Suppose that \textbf{\emph{(H$_\text{w}^\text{sub}$)}} holds. Let $\a\in(\ff{2}{q_b}+\ff{3\k-1}{1-\k},\a_b]$ and $q\in(\ff{2-2\k}{1+(1-\k)\a-3\k},q_b]$. Then for any $\de\in\mR$, $f\in\mL_T^q\bC_{\bba}^\a(\rho_\de)$ and $\l\geq0$, there exists a unique solution $u$ to PDE \eqref{weighted distributional PDE} such that for any $\t\in[0,2-\ff{2}{q})$,
\begin{align}\label{6.12}
(1+\l)^{\ff{\t}{2}}\|u\|_{\mL_T^\infty\bC_{\bba}^{2+\a-2/q-\t}(\rho_\nu)}\lesssim_C \|f\|_{\mL_T^q\bC_{\bba}^\a(\rho_\de)},
\end{align}
where $\nu:=\ff{2\k}{1+\a-2/q}+\de$ and $C:=C(\Xi,\de,\a,q,\t)>0$ is independent of $\l$.
\et

\bpf
We only prove the a priori estimate \eqref{6.12}. For any $r>0$ and $z_0\in\mR^{2d}$, we define
$$\phi_r^{z_0}(z):=\chi\Big(\ff{x-x_0}{(r(1+|z_0|_\bba))^3},\ff{v-v_0}{r(1+|z_0|_\bba)}\Big),$$
where $\chi\in C_c^\infty(\mR^{2d})$ with $\chi\equiv1$ on $B_{1/8}$, $\chi\equiv0$ on $B_{1/4}^c$, and 
$$(u^{z_0},b^{z_0},f^{z_0}):=(u\phi_r^{z_0},b\phi_{2r}^{z_0},f\phi_{2r}^{z_0}),\ \wt u^{z_0}:=u\phi_{2r}^{z_0}.$$
It is easy to see that 
$$\phi_{2r}^{z_0}(z)=1,\ \forall z\in\supp(\phi_r^{z_0}).$$
By \cite[Lemma 3.7]{HZZZ24}, for fixed small $r>0$, we have for any $\a\geq0$ and $\g\in\mR$,
\begin{align}\label{DZ1}
\|u(t)\|_{\bC_\bba^\a(\rho_\g)}\asymp\sup_{z_0\in\mR^{2d}}(\rho_\g(z_0)\|u^{z_0}(t)\|_{\bC_\bba^\a})\asymp\sup_{z_0\in\mR^{2d}}(\rho_\g(z_0)\|\wt u^{z_0}(t)\|_{\bC_\bba^\a}),
\end{align}
and for any $j\in\mN$,
\begin{align}\label{6.14-1}
\|\nn_v^j\phi_r^{z_0}\|_{\bC_\bba^\a}+\|v\cdot\nn_x\phi_r^{z_0}\|_{\bC_\bba^\a}\lesssim\rho(z_0),
\end{align}
\begin{align}\label{6.14-2}
\|\phi_r^{z_0}\rho_\g\|_{\bC_\bba^\a}\lesssim\rho_\g(z_0).
\end{align}
Hence,
\begin{align}\label{6.15-1}
\|b^{z_0}\|_{\mL_T^{q_b}\bC_\bba^{\a_b}}\overset{\eqref{paraproductcontrol4}}{\lesssim}\|b\rho_\k\|_{\mL_T^{q_b}\bC_\bba^{\a_b}}\|\phi_{2r}^{z_0}\rho_{-\k}\|_{\bC_\bba^1}\overset{\eqref{6.14-2}}{\lesssim}\rho_{-\k}(z_0)c_b,
\end{align}
\begin{align}\label{6.15-2}
\|f^{z_0}\|_{\mL_T^q\bC_\bba^\a}\overset{\eqref{paraproductcontrol4}}{\lesssim}\|f\rho_\de\|_{\mL_T^q\bC_\bba^\a}\|\phi_{2r}^{z_0}\rho_{-\de}\|_{\bC_\bba^1}\overset{\eqref{6.14-2}}{\lesssim}\rho_{-\de}(z_0)\|f\|_{\mL_T^q\bC_{\bba}^\a(\rho_\de)},
\end{align}
and for $\k<1$,
\begin{align}
&\|\div_v b^{z_0}\|_{\mL_T^{q_b}\bC_\bba^{\a_b}}\leq\|\phi_{2r}^{z_0}\div_v b\|_{\mL_T^{q_b}\bC_\bba^{\a_b}}+\|b\cdot\nn_v\phi_{2r}^{z_0}\|_{\mL_T^{q_b}\bC_\bba^{\a_b}}\no\\
&\qquad\overset{\eqref{paraproductcontrol4}}{\lesssim}\|\div_v b\|_{\mL_T^{q_b}\bC_\bba^{\a_b}}\|\phi_{2r}^{z_0}\|_{\bC_\bba^1}+\|b\rho_\k\|_{\mL_T^{q_b}\bC_\bba^{\a_b}}\|\phi_{4r}^{z_0}\rho_{-\k}\|_{\bC_\bba^1}\|\nn_v\phi_{2r}^{z_0}\|_{\bC_\bba^1}\no\\
&\quad\overset{\eqref{6.14-1},\eqref{6.14-2}}{\lesssim} c_b+c_b\rho_{1-\k}(z_0)\lesssim c_b.
\end{align}
Multiplying both sides of PDE \eqref{weighted distributional PDE} by $\phi_r^{z_0}$, we have
$$\p_t u^{z_0}=\Delta_v u^{z_0}-v\cdot\nn_x u^{z_0}-\l u^{z_0}+b^{z_0}\cdot\nn_v u^{z_0}+g^{z_0},$$
where
$$g^{z_0}=f^{z_0}-2\nn_v\wt u^{z_0}\cdot\nn_v \phi_r^{z_0}-\wt u^{z_0}\Delta_v\phi_r^{z_0}+v\cdot\nn_x\phi_r^{z_0} \wt u^{z_0}-b^{z_0}\cdot\nn_v\phi_r^{z_0}\wt u^{z_0}.$$
By \eqref{6.14-1}, \eqref{6.14-2}.\eqref{6.15-1} and \eqref{6.15-2}, for any $\a\leq\a_b$ and $q\leq q_b$, we have
\begin{align}
\|g^{z_0}\|_{\mL_T^q\bC_\bba^\a}&\ \leq\|f^{z_0}\|_{\mL_T^q\bC_\bba^\a}+\|2\nn_v\wt u^{z_0}\cdot\nn_v \phi_r^{z_0}-\wt u^{z_0}\Delta_v\phi_r^{z_0}\|_{\mL_T^q\bC_\bba^\a}\no\\
&\ \quad+\|v\cdot\nn_x\phi_r^{z_0} \wt u^{z_0}\|_{\mL_T^q\bC_\bba^\a}+\|b^{z_0}\cdot\nn_v\phi_r^{z_0}\wt u^{z_0}\|_{\mL_T^q\bC_\bba^{\a_b}}\no\\
&\overset{\eqref{paraproductcontrol4}}{\lesssim}\|f^{z_0}\|_{\mL_T^q\bC_\bba^\a}+\|\nn_v\phi_r^{z_0}\|_{\bC_\bba^1}\|\wt u^{z_0}\|_{\mL_T^q\bC_\bba^{\a+1}}+\|v\cdot\nn_x\phi_r^{z_0}\|_{\bC_\bba^1}\|\wt u^{z_0}\|_{\mL_T^q\bC_\bba^\a}\no\\
&\ \quad+\|b^{z_0}\|_{\mL_T^{q_b}\bC_\bba^{\a_b}}\|\nn_v\phi_r^{z_0}\|_{\bC_\bba^1}\|\wt u^{z_0}\|_{\mL_T^{q_0}\bC_\bba^1}\no\\
&\ \lesssim\rho_{-\de}(z_0)\|f\|_{\mL_T^q\bC_{\bba}^\a(\rho_\de)}+\rho(z_0)\|\wt u^{z_0}\|_{\mL_T^q\bC_\bba^{\a+1}}+\rho_{1-\k}(z_0)c_b\|\wt u^{z_0}\|_{\mL_T^{q_0}\bC_\bba^1}\no\\
&\ \lesssim\rho_{-\de}(z_0)\|f\|_{\mL_T^q\bC_{\bba}^a(\rho_\de)}+\rho_{1-\k}(z_0)\|\wt u^{z_0}\|_{\mL_T^{q_0}\bC_\bba^1},\no
\end{align}
where $q_0:=\ff{q_b q}{q_b-q}$. By \eqref{6.7}, for any $\t\in[0,2-\ff{2}{q})$ and $\l\geq0$, we have
\begin{align}
(1+\l)^{\ff{\t}{2}}\|u^{z_0}\|_{\mL_T^\infty\bC_\bba^{1+\b(\a,q)-\t}}&\lesssim(1+\|b^{z_0}\|_{\mL_T^{q_b}\bC_\bba^{\a_b}}+\|\div_v b^{z_0}\|_{\mL_T^{q_b}\bC_\bba^{\a_b}})^{\ff{2}{\b(\a,q)}}\|g^{z_0}\|_{\mL_T^q\bC_\bba^\a}\no\\
&\lesssim\rho_{-\ff{2\k}{\b(\a,q)}}(z_0)\Big(\rho_{-\de}(z_0)\|f\|_{\mL_T^q\bC_{\bba}^\a(\rho_\de)}+\rho_{1-\k}(z_0)\|\wt u^{z_0}\|_{\mL_T^{q_0}\bC_\bba^1}\Big)\no\\
&\lesssim\rho_{-\nu}(z_0)(\|f\|_{\mL_T^q\bC_{\bba}^\a(\rho_\de)}+\rho_{1-\k+\de}(z_0)\|\wt u^{z_0}\|_{\mL_T^{q_0}\bC_\bba^1})\no,
\end{align}
where $\b(a,q):=1+\a-\ff{2}{q}$. Multiplying both sides by $\rho_{\nu}(z_0)$ and 
taking supremum of $z_0$, by \eqref{DZ1}, we have for any $\l\geq0$,
\begin{align}\label{pf-3.29}
(1+\l)^{\ff{\t}{2}}\|u\|_{\mL_T^\infty\bC_{\bba}^{1+\b(\a,q)-\t}(\rho_\nu)}\lesssim\|f\|_{\mL_T^q\bC_{\bba}^\a(\rho_\de)}+\|u\|_{\mL_T^{q_0}\bC_{\bba}^1(\rho_{1-\k+\de})}.
\end{align}
It follows from \eqref{interpolation} and Young's inequality that for any small $\varepsilon>0$,
\begin{align}
\|u\|_{\mL_T^{q_0}\bC_{\bba}^1(\rho_{1-\k+\de})}\lesssim\varepsilon\|u\|_{\mL_T^{q_0}\bC_{\bba}^{1+\b(\a,q)-\t}(\rho_{1-\k+\de})}+C_\varepsilon\|u\|_{\mL_T^{q_0}\bC_{\bba}^0(\rho_{1-\k+\de})}.\no
\end{align}
Substituting this back into \eqref{pf-3.29} and noting $1-\k+\de\geq\nu$, we have
\begin{align}\label{6.18}
(1+\l)^{\ff{\t}{2}}\|u\|_{\mL_T^\infty\bC_{\bba}^{1+\b(\a,q)-\t}(\rho_\nu)}\lesssim\|f\|_{\mL_T^q\bC_{\bba}^\a(\rho_\de)}+\|u\|_{\mL_T^{q_0}\bC_{\bba}^0(\rho_{1-\k+\de})}.
\end{align}
On the other hand, for any $\l\geq0$ we have
\begin{align}
\|u^{z_0}\|_{\mL_T^\infty\bC_\bba^0}&\lesssim\|u^{z_0}\|_{\mL_T^\infty}\overset{\eqref{6.8}}{\lesssim}\rho_{\ff{\b(\a,q)-1}{\b(\a,q)}\k}(z_0)\|g^{z_0}\|_{\mL_T^q\bC_\bba^\a}\no\\
&\lesssim\rho_{\ff{\b(\a,q)-1}{\b(\a,q)}\k-\de}(z_0)\|f\|_{\mL_T^q\bC_{\bba}^\a(\rho_\de)}+\rho_{1-\ff{\k}{\b(\a,q)}}(z_0)\|\wt u^{z_0}\|_{\mL_T^{q_0}\bC_\bba^1}\no\\
&\lesssim\rho_{-1+\k-\de}(z_0)(\rho_{1-\ff{\k}{\b(\a,q)}}(z_0)\|f\|_{\mL_T^q\bC_{\bba}^\a(\rho_\de)}+\rho_{2-\k-\ff{3\k}{\b(\a,q)}}(z_0)\|u\|_{\mL_T^{q_0}\bC_{\bba}^1(\rho_\nu)}).\no
\end{align}
Note that $2-\k-\ff{3\k}{\b(\a,q)}>0$. As above, by \eqref{DZ1}, we obtain
\begin{align}\label{3.31-17}
\|u\|_{\mL_T^\infty\bC_{\bba}^0(\rho_{1-\k+\de})}\lesssim\|f\|_{\mL_T^q\bC_{\bba}^\a(\rho_\de)}+\|u\|_{\mL_T^{q_0}\bC_{\bba}^1(\rho_\nu)}\overset{\eqref{6.18}}{\lesssim}\|f\|_{\mL_T^q\bC_{\bba}^\a(\rho_\de)}+\|u\|_{\mL_T^{q_0}\bC_{\bba}^0(\rho_{1-\k+\de})}.
\end{align}
For any $t\in[0,T]$ and $q\in(\ff{2-2\k}{1+(1-\k)\a-3\k},q_b)$, we have 
\begin{align}\label{3.29-17}
\|u(t)\|^{q_0}_{\bC_{\bba}^0(\rho_{1-\k+\de})}\lesssim\int_0^t\|u(s)\|^{q_0}_{\bC_{\bba}^0(\rho_{1-\k+\de})}\dif s+\|f\|^{q_0}_{\mL_T^q\bC_{\bba}^\a(\rho_\de)}.
\end{align}
It follows from Gronwall's inequality that
\begin{align}\label{3.30-17}
\|u\|_{\mL_T^\infty\bC_{\bba}^0(\rho_{1-\k+\de})}\lesssim\|f\|_{\mL_T^q\bC_{\bba}^\a(\rho_\de)}.
\end{align}
Plugging \eqref{3.30-17} back into \eqref{3.31-17}, we have \eqref{3.30-17} still holds for $q=q_b$. Finally, substituting \eqref{3.30-17} into \eqref{6.18}, we complete the proof.
\epf

{\section{Solvability of SDEs with distributional drifts}\label{sec4}}

Now that we have established the well-posedness of the associated PDE \eqref{weighted distributional PDE} with distributional drifts, we turn our attention to the weak well-posedness of the SDE \eqref{3} under some subcritical assumption. Throughout this section, we fix $T>0$ and define $\Gamma_n(z)=\Gamma_n(x,v):=n^{4d}\Gamma(nx,n^3v)$, where $\Gamma\in C_c^\infty(\mR^{2d})$ is a smooth probability density function. In the subsequent analysis, unless otherwise stated, we always assume that
\begin{enumerate}[\textbf{(H$^\text{sub}$)}]
\item $b=b_1+b_2$, where the singular component $b_1$ satisfies that
\begin{align}\label{ellb1}
\ell_{b_1}:=\|b_1\|_{\mL_T^{q_b}\bC_\bba^{\a_b}}+\|\div_v b_1\|_{\mL_T^{q_b}\bC_\bba^{\a_b}}<\infty,
\end{align}
with $\a_b\in(-1,0)$, $q_b\in(\ff{2}{1+\a_b},\infty]$, and the regular component $b_2$ satisfies the linear growth condition:  for some $c_0,c_1\geq0$ and $l\in L_T^2:=L^2([0,T];\mR_+)$,
\begin{align}\label{lineargrowth}
|b_2(t,z)|\leq l(t)(c_0+c_1|z|_\bba).
\end{align}
\end{enumerate}

Denote the parameter set:
$$\Xi':=(T,d,\a_b,q_b,\ell_{b_1},c_0,c_1,\| l\|_{L_T^2}).$$

\br\label{GubinelliDecomposition}
By a similar argument of \cite[Appendix. A]{HZ23}, for any $b\in\bC_\bba^{\a_b}(\rho_\k)$ with $\k\in[0,1+\a_b)$, we have the decomposition $b=b_1+b_2$, where $b_1\in\bC_\bba^{\a_b-\k}$ belongs to unweighted anisotropic H\"{o}lder space and $b_2\in\bC_\bba^{\a_b+1-\k}(\rho_1)$ satisfies the linear growth condition. Thus, the assumption \textbf{\emph{(H$^\text{sub}$)}} is weaker than assumption {\textbf{\emph{(H$^\text{sub}_{\bf w}$)}}}.
\er

\subsection{Krylov's estimate and Young's integral}\label{sec41}

We first introduce the following crucial Krylov's estimate.

\bd
Let $(Z_t)_{t\in[0,T]}$ be a $\mR^{2d}$-valued stochastic process on a probability space $(\Omega,\sF,\bP)$ and consider the space $\mL_T^q\bC_\bba^\a(\rho_\de)$ with parameters $q\in(1,\infty]$, $\a\in(-1,0)$ and $\de\in\mR$. For given $\t\in(0,1)$ and $p\geq2$
, one says that the pair $(Z,\mL_T^q\bC_\bba^\a(\rho_\de))$ satisfies the Krylov's estimate with parameters $(p,\t)$ and constant $C>0$ if for any $f\in\mL_T^1 C_b(\rho_\de)\cap\mL_T^q\bC_\bba^\a(\rho_\de)$ and stopping times $0\leq \tau_0\leq \tau_1\leq T$ with $\tau_1-\tau_0\leq\s$,
\begin{align}\label{krylov}
\Big\|\int_{\tau_0}^{\tau_1}f(s,Z_s)\dif s\Big\|_{L^p(\Omega)}\lesssim_C\s^{\ff{1+\t}{2}}\|f\|_{\mL_T^q\bC_\bba^\a(\rho_\de)}.
\end{align}
\ed

As a direct consequence of Krylov's estimate and Kolmogorov continuity theorem, we derive the following result (cf. \cite[Proposition 4.3]{HZ23}).

\bp\label{limit}
For given $\t\in(0,1)$ and $p\geq2$, suppose that $(Z,\mL_T^q\bC_{\bba}^\a(\rho_\de))$ satisfies the Krylov's estimate with parameter $(p,\t)$. For any $f\in\mL_T^q\bC_{\bba}^\a(\rho_\de)$, if there exist $f_n\in \mL_T^1C_b(\rho_\de)\cap\mL_T^q\bC_{\bba}^\a(\rho_\de)$ such that
$$\lim_{n\to\infty}\|f_n-f\|_{\mL_T^q\bC_{\bba}^\a(\rho_\de)}=0,$$
then the limit $\lim_{n\to\infty}\int_0^.f_n(s,Z_s)\dif s$ exists in $L^p(\Omega,C([0,T]))$ and we denote the limit by $A_t^f:=\lim_{n\to\infty}\int_0^t f_n(s,Z_s)\dif s$, which is independent of the choice of $\{f_n\}_{n\in\mN}$. Furthermore, we have the estimate
\begin{align}\label{4.1}
\sup_{t\in[0,T]}\|A_t^f\|_{L^p(\Omega)}+\sup_{s\neq t\in[0,T]}\ff{\|A_t^f-A_s^f\|_{L^p(\Omega)}}{|t-s|^{(1+\t)/2}}\lesssim\|f\|_{\mL_T^q\bC_{\bba}^\a(\rho_\de)}.
\end{align}
\ep

\bpf
Note that $\{f_n\}_{n\in\mN}$ is a Cauchy sequence in $\mL_T^q\bC_{\bba}^\a(\rho_\de)$. Since $\ff{(1+\t)p}{2}>1$, by Krylov's estimate and Kolmogorov continuity theorem, we conclude that $\{\int_0^.f_n(s,Z_s)\dif s\}_{n\in\mN}$ is also a Cauchy sequence in $L^p(\Omega,C([0,T]))$, which implies that the limit $\lim_{n\to\infty}\int_0^.f_n(s,Z_s)\dif s$ exists. Furthermore, \eqref{4.1} follows from \eqref{krylov} directly.
\epf

We also need the following proposition about the substitution for Young's integrals (see \cite[Proposition 4.4]{HZ23}).

\bp\label{pp4.4}
Let $q\in(1,\infty]$, $\a\in(-1,0)$ and $\de\in\mR$ be given parameters. For $p>2$ and $\theta\in(0,1)$ with $\theta p>2$, suppose that $(Z,\mL_T^q\bC_{\bba}^\a(\rho_\de))$ satisfies the Krylov's estimate with parameter $(p,\t)$ and for any $\b\in(0,\ff{1}{2})$,
\begin{align}
\bigg\|\sup_{s\neq t\in[0,T]}\ff{|Z_t-Z_s|}{|t-s|^\b}\bigg\|_{L^p(\Omega)}<\infty.
\end{align}
Let $g:\mR_+\times\mR^{2d}\to\mR$ be a bounded function satisfying that
\begin{align}\label{4.3}
|g(t,z_1)-g(s,z_2)|\lesssim_C\sqrt{|t-s|}+|z_1-z_2|.
\end{align}
Then for any $f\in\cup_{\varepsilon>0}\mL_T^q\bC_{\bba}^{\a+\varepsilon}(\rho_\de)$, the integral $\int_0^tg(s,Z_s)\dif A_s^f$ is well-defined as Young's integral and 
\begin{align}\label{4.4}
\int_0^tg(s,Z_s)\dif A_s^f=A_t^{f\cdot g},\ a.s.,
\end{align}
where $A_\cdot^f$ is defined in Proposition \ref{limit}.
\ep

Below we consider the following approximation SDE:
\begin{align}\label{approximation}
\dif Z_t^n=B_n(t,Z_t^n)\dif t+\s\dif W_t,\ Z_0^n=z_0=(x_0,v_0),
\end{align}
where the drift term is given by $B_n(t,x,v):=(v,(b_{1;n}+b_2)(t,x,v))^{\rm t}$ with $b_{1;n}:=b_1*\G_n$ denoting the standard mollification approximation. By \eqref{molifycompare} we have $\ell_{b_{1;n}}\leq\ell_{b_1}$.

Now we can show the following crucial Krylov's estimate for the solutions of the approximation SDE \eqref{approximation}.

\bl\label{kryloves}
Suppose that $Z^n$ is the unique weak solution of SDE \eqref{approximation}. Then we have
\begin{itemize}
\item[(i)] Under {\bf(H$^\text{sub}$)}, for any $\a\in(\ff{2}{q_b}-1,\a_b]$, $q\in(\ff{2}{1+\a},q_b]$ and $p\geq 2$, $(Z^n,\mL_T^q\bC_\bba^\a)$ satisfies the Krylov's estimate \eqref{krylov} with parameters $(p,1+\a-\ff{2}{q})$ and constant $C:=C(\Xi',\a,q,p)>0$.
\item[(ii)] Under \textbf{\emph{(H$_\text{w}^\text{sub}$)}}, for any $\de\in\mR$, $\a\in(\ff{2}{q_b}+\ff{3\k-1}{1-\k},\a_b]$, $q\in(\ff{2-2\k}{1+(1-\k)\a-3\k},q_b]$ and $p\geq 2$, $(Z^n,\mL_T^q\bC_\bba^\a(\rho_\de))$ satisfies the Krylov's estimate \eqref{krylov} with parameters $(p,1+\a-\ff{2}{q})$ and constant $$C:=C_0(1+|z_0|_\bba)^{\ff{2\k}{1+\a-2/q}+\de},$$
where $C_0:=C_0(\Xi,\de,\a,q,p)>0$.
\end{itemize}
\el

\br
In the second part of Lemma \ref{kryloves}, we strengthen our hypotheses by replacing the assumption {\bf(H$^\text{sub}$)} with stronger assumption {\bf(H$_\text{w}^\text{sub}$)} to establish the Krylov's estimate on weighted spaces, which enables us to consider the asymptotic behavior - either growth or decay properties - of the probability density associated with weak solutions to SDE \eqref{3}.
\er

To prove the above Krylov's estimate, we need the following uniform moment estimate.

\bl
Under {\bf(H$^\text{sub}$)}, for any $z_0\in\mR^{2d}$ and $\de\in\mR$, there exists a constant $C:=C(\Xi',\de)>0$ such that
\begin{align}\label{unimomentes}
\sup_{n}\bE\bigg(\sup_{t\in[0,T]}\rho_\de(Z_t^n)\bigg)\lesssim\rho_\de(z_0).
\end{align}
\el

\bpf
Let $u_n\in \mL^\infty_T C^2_b$ solve the following backward PDE  (cf. \cite[Theorem 1.3]{CHM21}):
\begin{align}
\p_t u_n+\Delta_v u_n+v\cdot\nn_x u_n+b_{1;n}\cdot\nn_x u_n+b_{1;n}=0,\ u_n(T)=0.\no
\end{align}
By \eqref{6.12} we have 
\begin{align}\label{KG1}
\|\nn_v u_n\|_\infty\vee\|u_n\|_\infty\lesssim\|u_n\|_{\mL_T^\infty\bC_\bba^1}\lesssim\|b_{1;n}\|_{\mL_T^{q_b}\bC_\bba^{a_b}}\lesssim\|b_1\|_{\mL_T^{q_b}\bC_\bba^{a_b}}\leq\ell_{b_1}. 
\end{align}
Here and below, all the implicit constants are independent of $n$. For $Z_t^n=(X_t^n,V_t^n)$, define $\wt Z_t^n:=(X_t^n,V_t^n+u_n(t,X_t^n,V_t^n))$. By It\^{o}'s formula, we have 
\begin{align}
\dif\wt Z_t^n=\wt B_n(t,Z_t^n)\dif t+\wt\s_n(t,Z_t^n)\dif W_t,\ \wt Z_0^n=\wt z_0=(x_0,v_0+u_n(0,x_0,v_0)),
\end{align}
with
\begin{align}
\wt B_n(s,z):=\wt B_n(s,x,v)=(v,b_2\cdot(\mI+\nn_v u_n)(s,z))^{\rm t},\no
\end{align}
\begin{align}
\wt \s_n(s,z):=(0,\sqrt{2}(\mI+\nn_v u_n)(s,z))^{\rm t}.\no
\end{align}
Moreover, it holds that
\begin{align}\label{wtbn}
|\wt B_n(s,z)|\lesssim l(s)(c_0+c_1|z|_\bba),
\end{align}
and there exists a constant $\k_0>0$ such that for any $\xi\in\mR^{2d}$,
\begin{align}\label{wtsigman}
\k_0^{-1}|\xi|^2\leq|\wt\s_n(s,z)\xi|^2\leq\k_0|\xi|^2.
\end{align}
Note that by \eqref{KG1},
\begin{align}
1+|\wt Z_t^n|_\bba\lesssim1+|Z_t^n|_\bba+\|u_n\|_\infty\lesssim1+|Z_t^n|_\bba\lesssim1+|\wt Z_t^n|_\bba+\|u_n\|_\infty\lesssim1+|\wt Z_t^n|_\bba,\no
\end{align}
which implies that for any $\de\in\mR$, 
$$
\rho_\de(\wt Z_t^n)\asymp\rho_\de(Z_t^n).
$$
Hence, it suffices to prove that 
\begin{align}\label{wtznuni}
\sup_{n}\bE\bigg(\sup_{t\in[0,T]}\rho_\de(\wt Z_t^n)\bigg)\lesssim\rho_\de(\wt z_0).
\end{align}
Applying It\^{o}'s formula, we have for any $n\in\mN$,
\begin{align}
\rho_\de(\wt Z_t^n)=\rho_\de(\wt z_0)+\int_0^t(\tr(\wt\s_n\wt\s_n^{\rm t})(s,Z_s^n)\cdot\nn^2+\wt B_n(s,Z_s^n)\cdot\nn)\rho_\de(\wt Z_s^n)\dif s+\int_0^t(\wt \s_n(s,Z_s^n)\cdot\nn)\rho_\de(\wt Z^n_s)\dif W_s.\no
\end{align}
By BDG's inequality, we have
\begin{align}
\bE\bigg(\sup_{t\in[0,T]}\rho_{2\de}(\wt Z_t^n)\bigg)&=\bE\bigg(\sup_{t\in[0,T]}(\rho_{\de}(\wt Z_t^n))^2\bigg)\no\\
&\lesssim\rho_{2\de}(\wt z_0)+\bE\Big(\int_0^T|(\tr(\wt\s_n\wt\s_n^{\rm t})(s,Z_s^n)\cdot\nn^2+\wt B_n(s,Z_s^n)\cdot\nn)\rho_{\de}(\wt Z_s^n)|\dif s\Big)^2\no\\
&\quad+\bE\Big(\int_0^T|\wt\s_n(s,Z_s^n)\cdot\nn\rho_\de(\wt Z_s^n)|^2\dif s\Big).
\end{align}
By \eqref{1.14} we have
\begin{align}
|\tr(\wt\s_n\wt\s_n^{\rm t})(s,Z_s^n)\cdot\nn^2)\rho_{\de}(\wt Z_s^n)|\overset{\eqref{wtsigman}}{\lesssim}|\nn^2\rho_\de(\wt Z_s^n)|\lesssim\rho_{\de+2}(\wt Z_s^n)\lesssim\rho_\de(\wt Z_s^n),\no
\end{align}
\begin{align}
|(\wt B_n(s,Z_s^n)\cdot\nn)\rho_{\de}(\wt Z_s^n)|\lesssim|\wt B_n(s,Z_s^n)||\nn\rho_\de(\wt Z_s^n)|\overset{\eqref{wtbn}}{\lesssim}l(s)\rho_{-1}(Z_s^n)\rho_{\de+1}(\wt Z_s^n)\lesssim l(s)\rho_\de(\wt Z_s^n),\no
\end{align}
\begin{align}
|\wt\s_n(s,Z_s^n)\cdot\nn\rho_\de(\wt Z_s^n)|^2\overset{\eqref{wtsigman}}{\lesssim}|\nn\rho_\de(\wt Z_s^n)|^2\lesssim\rho_{2\de+2}(\wt Z_s^n)\lesssim\rho_{2\de}(\wt Z_s^n).\no
\end{align}
Hence,
\begin{align}
\bE\bigg(\sup_{t\in[0,T]}\rho_{2\de}(\wt Z_t^n)\bigg)\lesssim\rho_{2\de}(\wt z_0)+\int_0^T(1+l(s))^2\bE\rho_{2\de}(\wt Z_s^n)\dif s.
\end{align}
Since $\de\in\mR$ is arbitrary, by Gronwall's inequality, we obtain \eqref{wtznuni} and complete the proof.
\epf

Now we come to the proof of Lemma \ref{kryloves}.

\bpf[Proof of Lemma \ref{kryloves}]
(i) For $f\in C_c^\infty([0,T]\times\mR^{2d})$, let $u_n$ solve the following backward PDE:
\begin{align}
\p_t u_n+\Delta_v u_n+v\cdot\nn_x u_n-\l u_n+b_{1;n}\cdot\nn_v u_n=f,\ w_n(T)=0.\no
\end{align}
By It\^{o}'s formula, we have
\begin{align*}
u_n(t,Z_t^n)=u_n(0,z_0)+\int_0^t(\l u_n-f)(s,Z_s^n)\dif s+\int_0^t b_2\cdot\nn_v u_n(s,Z_s^n)\dif s+\sqrt{2}\int_0^t \nn_v u_n(s,Z_s^n)\dif W_s.
\end{align*}
For any $p\geq2$ and stopping times $0\leq \tau_0\leq \tau_1\leq T$ with $\tau_1-\tau_0\leq\de$, by BDG's inequality we have
\begin{align}
\Big\|\int_{\tau_0}^{\tau_1} f(s,Z_s^n)\dif s\Big\|_{L^p(\Omega)}&\lesssim\|u_n(\tau_1,Z_{\tau_1}^n)-u_n(\tau_0,Z_{\tau_0}^n)\|_{L^p(\Omega)}+\Big\|\int_{\tau_0}^{\tau_1}\l u_n(s,Z_s^n)\dif s\Big\|_{L^p(\Omega)}\no\\
&\quad+\Big\|\Big(\int_{\tau_0}^{\tau_1}|\nn_v u_n(s,Z_s^n)|^2\dif s\Big)^{\ff{1}{2}}\Big\|_{L^p(\Omega)}\no\\
&\lesssim2\|u_n\|_\infty+\l\s\|u_n\|_\infty+\|\nn_v u_n\|_\infty\Big\|\int_{\tau_0}^{\tau_1}l(s)(c_0+c_1|Z_s^n|_\bba)\dif s\Big\|_{L^p(\Omega)}\no\\
&\quad+\Big\|\Big(\int_{\tau_0}^{\tau_1}|\nn_v u_n(s,Z_s^n)|^2\dif s\Big)^{\ff{1}{2}}\Big\|_{L^p(\Omega)}\no\\
&\lesssim(1+\l\s)\|u_n\|_\infty+\|\nn_v u_n\|_\infty\s^{\ff{1}{2}}\| l\|_{L_T^2}\Big\| \sup_{t\in[0,T]}(c_0+c_1|Z_t^n|_\bba)\Big\|_{L^p(\Omega)}\no\\
&\quad+\s^{\ff{1}{2}}\|\nn_v u_n\|_\infty.\no
\end{align}
By \eqref{unimomentes}, we have
\begin{align}
\Big\|\int_{\tau_0}^{\tau_1}f(s,Z_s^n)\dif s\Big\|_{L^p(\Omega)}\lesssim(1+\l\s)\|u_n\|_\infty+(\s(1+|z_0|_\bba)+\s^{\ff{1}{2}})\|\nn_v u_n\|_\infty.\no
\end{align}
By \eqref{6.12} with unweighted Schauder's estimates, we have
\begin{align}
\|u_n\|_\infty\lesssim\|u_n\|_{\mL_T^\infty\bC_\bba^0}\lesssim(1+\l)^{-\ff{2+\a-2/q}{2}}\|f\|_{\mL_T^q\bC_\bba^\a},\no
\end{align}
\begin{align}
\|\nn_v u_n\|_\infty\lesssim\|u_n\|_{\mL_T^\infty\bC_\bba^1}\lesssim(1+\l)^{-\ff{1+\a-2/q}{2}}\|f\|_{\mL_T^q\bC_\bba^\a}.\no
\end{align}
Thus,
\begin{align}
\Big\|\int_{\tau_0}^{\tau_1}f(s,Z_s^n)\dif s\Big\|_{L^p(\Omega)}&\lesssim\Big[(1+\l\s)(1+\l)^{-\ff{2+\a-2/q}{2}}+\s^{\ff{1}{2}}(1+\l)^{-\ff{1+\a-2/q}{2}}\Big]\|f\|_{\mL_T^q\bC_\bba^\a}\no\\
&\quad+\s(1+|z_0|_\bba)(1+\l)^{-\ff{1+\a-2/q}{2}}\|f\|_{\mL_T^q\bC_\bba^\a}\no.
\end{align}
For any $\s\in(0,1)$, taking $\l=\s^{-1}$ we get
\begin{align}
\Big\|\int_{\tau_0}^{\tau_1}f(s,Z_s^n)\dif s\Big\|_{L^p(\Omega)}\lesssim\s^{\ff{2+\a-2/q}{2}}\|f\|_{\mL_T^q\bC_\bba^\a}.\no
\end{align}
By a standard approximation, the above estimate still holds for $f\in\mL_T^1 C_b\cap\mL_T^q\bC_\bba^\a$.

(ii) For fixed $f\in C_c^\infty([0,T]\times\mR^{2d})$, let $u_n$ solve the following backward PDE:
$$\p_t u_n+\Delta_v u_n+v\cdot\nn_x u_n-\l u_n+b_n\cdot\nn_v u_n +f=0,\ u_n(T)=0.$$
Let $\nu=\ff{2\k}{1+\a-2/q}+\de$. As above, for stopping times $0\leq \tau_0\leq \tau_1\leq T$ with $\tau_1-\tau_0\leq\de$ we have
\begin{align}
\Big\|\int_{\tau_0}^{\tau_1} f(s,Z_s^n)\dif s\Big\|_{L^p(\Omega)}&\lesssim\|u_n(\tau_1,Z_{\tau_1}^n)-u_n(\tau_0,Z_{\tau_0}^n)\|_{L^p(\Omega)}+\Big\|\int_{\tau_0}^{\tau_1}\l u_n(s,Z_s^n)\dif s\Big\|_{L^p(\Omega)}\no\\
&\quad+\Big\|\Big(\int_{\tau_0}^{\tau_1}|\nn_v u_n(s,Z_s^n)|^2\dif s\Big)^{\ff{1}{2}}\Big\|_{L^p(\Omega)}=:I_1+I_2+I_3.\no
\end{align}
Note that
\begin{align}
I_1&\ \leq\|(u_n\rho_\nu)(\tau_1,Z_{\tau_1}^n)\rho_{-\nu}(Z_{\tau_1}^n)\|_{L^p(\Omega)}+\|(u_n\rho_\nu)(\tau_0,Z_{\tau_0}^n)\rho_{-\nu}(Z_{\tau_0}^n)\|_{L^p(\Omega)}\no\\
&\ \leq2\|u_n\rho_\nu\|_\infty\Big\|\sup_{t\in[0,T]}\rho_{-\nu}(Z_t^n)\Big\|_{L^p(\Omega)}\overset{\eqref{0117:01}}{\lesssim}\|u_n\|_{\mL_T^\infty\bC_{\bba}^0(\rho_\nu)}\Big[\bE\sup_{t\in[0,T]}\rho_{-p\nu}(Z_t^n)\Big]^{\ff{1}{p}}\no\\
&\overset{\eqref{unimomentes}}{\lesssim}(1+|z_0|_\bba)^{\nu}\|u_n\|_{\mL_T^\infty\bC_{\bba}^0(\rho_\nu)}\no.
\end{align}
Similarly, for $\l:=\s^{-1}$,
\begin{align}
I_2\lesssim\lambda\sigma\|u_n\rho_\nu\|_\infty\Big\|\sup_{t\in[0,T]}\rho_{-\nu}(Z_t^n)\Big\|_{L^p(\Omega)}\lesssim(1+|z_0|_\bba)^{\nu}\|u_n\|_{\mL_T^\infty\bC_{\bba}^0(\rho_\nu)}\no,
\end{align}
\begin{align}
I_3\lesssim \sigma^{\ff{1}{2}}\|\nn_v u_n\rho_\nu\|_\infty\Big\|\sup_{t\in[0,T]}\rho_{-\nu}(Z_t^n)\Big\|_{L^p(\Omega)}\lesssim(1+|z_0|_\bba)^{\nu}\s^{\ff{1}{2}}\|\nn_v u_n\|_{\mL_T^\infty\bC_{\bba}^0(\rho_\nu)}.\no
\end{align}
Hence,
\begin{align}
\Big\|\int_{\tau_0}^{\tau_1} f(s,Z_s^n)\dif s\Big\|_{L^p(\Omega)}&\quad\lesssim(1+|z_0|_\bba)^{\nu}\Big[\|u_n\|_{\mL_T^\infty\bC_{\bba}^0(\rho_\nu)}+\s^{\ff{1}{2}}\|\nn_v u_n\|_{\mL_T^\infty\bC_{\bba}^0(\rho_\nu)}\Big]\no\\
&\ \overset{\eqref{6.12}}{\lesssim}(1+|z_0|_\bba)^{\nu}\|f\|_{\mL_T^q\bC_{\bba}^\a(\rho_\de)}[(1+\l)^{-\ff{2+\a-2/q}{2}}+\s^{\ff{1}{2}}(1+\l)^{-\ff{1+\a-2/q}{2}}]\no\\
&\quad\lesssim(1+|z_0|_\bba)^{\nu}\|f\|_{\mL_T^q\bC_{\bba}^\a(\rho_\de)}\s^{\ff{2+\a-2/q}{2}}.\no
\end{align}
By a standard approximation, we complete the proof.
\epf

\subsection{Existence of the solutions to subcritical kinetic SDEs}\label{sec42}

Let $\mC_T:=C([0,T];\mR^{2d})$ be the Banach space of all continuous functions and denote a path in $\mC_T$ by $\omega$. The canonical process is denoted by $\omega_t$. Let $\sB^s_t$ denote the natural filtration generated by $\{\om_r;r\in[s,t]\}$ and define
\begin{align}
\sB_t:=\sB^0_t,\ \sB^s:=\sB^s_T.\no
\end{align}
In this subsection we establish the existence of weak solution in the sense of Definition \ref{weaksol} by using a standard weak convergence method. To start with, we present the following tightness result of the law $\mP_n$ of approximation solutions $Z^n$ defined in \eqref{approximation} in $\mC_T$.

\bl\label{tight}
Under \textbf{\emph{(H$^\text{sub}$)}}, the family of probability measures $(\mP_n)_{n\in\mN}$ is tight.
\el

\bpf
Let $\tau\geq0$ be any stopping time less than T. Note that by \eqref{approximation} we have for any  $t>0$,  
\begin{align}
Z_{\tau+t}^n-Z_\tau^n=\int_\tau^{\tau+t} B_n(s,Z_s^n)\dif s+\s(W_{\tau+t}-W_\tau),\no
\end{align}
which implies that
\begin{align}
|Z_{\tau+t}^n-Z_\tau^n|\leq\int_\tau^{\tau+t} |V_s^n|\dif s+\int_\tau^{\tau+t} |b_{1;n}(s,Z_s^n)|\dif s+\int_\tau^{\tau+t} |b_2(s,Z_s^n)|\dif s+|\s(W_{\tau+t}-W_\tau)|.\no
\end{align}
By Lemma \ref{kryloves}, \eqref{unimomentes} and BDG's inequality, we have for any $\de>0$,
\begin{align}
\bE\bigg(\sup_{t\in[0,\de]}|Z_{\tau+t}^n-Z_\tau^n|\bigg)&\leq \bE\int_\tau^{\tau+\de} \Big(|V_s^n|+l(s)(c_0+c_1|Z_s^n|_\bba)\Big)\dif s+\bE\int_\tau^{\tau+\de}|b_{1;n}|(s,Z_s^n)\dif s\no\\
&\quad+\sqrt{2}\bE\bigg(\sup_{t\in[0,\de]}|W_{\tau+t}-W_\tau|\bigg)\no\\
&\lesssim_C \de^{\ff{1}{2}}\bigg(1+\bE\sup_{t\in[0,T]}|Z_t^n|_\bba\bigg)+\de^{\ff{2+\a_b-2/q_b}{2}}\|b_{1;n}\|_{\mL_T^{q_b}\bC_{\bba}^{\a_b}}+\de^{\ff{1}{2}}\no\\
&\lesssim_C \de^{\ff{1}{2}}(1+|z_0|_\bba)+\de^{\ff{2+\a_b-2/q_b}{2}}\|b_{1;n}\|_{\mL_T^{q_b}\bC_{\bba}^{\a_b}},\no
\end{align}
where the constant $C$ is independent of $n$. Thus, by \cite[Lemma 2.7]{ZZ18}, we obtain
\begin{align}
\sup_n\bE\Bigg(\sup_{\substack{0\leq t<t'<t+\delta\leq T}}|Z_{t'}^n-Z_t^n|^{\ff{1}{2}}\Bigg)\lesssim\de^{\ff{1}{4}}(1+|z_0|_\bba)^{\ff{1}{2}}+\de^{\ff{2+\a_b-2/q_b}{4}}\|b_{1;n}\|_{\mL_T^{q_b}\bC_{\bba}^{\a_b}}^{\ff{1}{2}}.\no
\end{align}
Then by Chebyshev's inequality, we have for any $\varepsilon>0$,
\begin{align}
\lim_{\de\to0}\sup_n\bP\Bigg(\sup_{\substack{0\leq t<t'<t+\delta\leq T}}|Z_{t'}^n-Z_t^n|>\varepsilon\Bigg)=0.\no
\end{align}
Hence, by \cite[Theorem 1.3.2]{SV06}, we derive the tightness.
\epf

As an immediate consequence of Lemma \ref{tight} and the Prohorov's theorem (cf. \cite[p. 9, Theorem 5.4]{RY99}), the sequence $(\mP_n)_{n\in\mN}$ is relatively compact in $\sP(\mC_T)$. 
Next, we use the accumulation point of the sequence $(\mP_n)_{n\in\mN}$ to construct a weak solution to SDE \eqref{3}.

\bt\label{existence}
Under \textbf{\emph{(H$^\text{sub}$)}}, there exists a weak solution $(\mathfrak F,Z,W)$ to SDE \eqref{3} in the sense of Definition \ref{weaksol} such that for all $t\in[0,T]$ and $f\in\mL_T^1 C_b\cap\mL_T^q\bC_\bba^\a$,
\begin{align}\label{exies}
\Big\|\int_{t_0}^{t_1}f(s,Z_s)\dif s\Big\|_{L^p(\Omega)}\lesssim_C (t_1-t_0)^{\ff{2+\a-2/q}{2}}\|f\|_{\mL_T^q\bC_\bba^\a},
\end{align}
where the constant $C=C(\Xi',\a,q)>0$. Moreover, under  {\bf{(H$_\text{w}^\text{sub}$)}},  we also have for any $\de\in\mR$,
\begin{align}\label{exies2}
\Big\|\int_{t_0}^{t_1}f(s,Z_s)\dif s\Big\|_{L^p(\Omega)}\lesssim_C \rho_{-\ff{2\k}{1+\a-2/q}-\de}(z_0)(t_1-t_0)^{\ff{2+\a-2/q}{2}}\|f\|_{\mL_T^q\bC_{\bba}^\a(\rho_\de)}.
\end{align}
\et

\bpf
Denote $\wt{\mathfrak F}:=(\wt\Omega,\wt\sF,(\wt\sF_t),\wt\bP)$ and let $(\wt{\mathfrak F},\wt Z^n,\wt W)$ be the weak solution to SDE \eqref{approximation}. By Lemma \ref{tight} and Prohorov's theorem, there is a subsequence still denoted by $n$ such that $\mQ_n:=\wt\bP\circ(\wt Z^n_\cdot,\wt W_\cdot)^{-1}\in\sP(\mC_T\times\mC_T)$ weakly converges to some probability measure $\mQ$. By Skorokhod's representation theorem, there is a probability space $\mathfrak F:=(\Omega,\sF,(\sF_t),\bP)$ and $\mC_T\times\mC_T$-valued random variables $(Z^n,W^n)$ and $(Z,W)$ defined on it such that 
$$
\bP\circ(Z^n,W^n)^{-1}=\mQ_n,\ \bP\circ(Z,W)^{-1}=\mQ
$$ 
and
\begin{align}
(Z^n,W^n)\rightarrow(Z,W)\ \bP-\text{a.s.},
\end{align}
\begin{align}
\dif Z_t^n=B_n(t,Z_t^n)\dif t+\s\dif W_t^n.
\end{align}
By (i) of Lemma \ref{kryloves}, we have
\begin{align}
\Big\|\int_{t_0}^{t_1}f(s,Z_s^n)\dif s\Big\|_{L^p(\Omega)}\lesssim_C (t_1-t_0)^{\ff{2+\a-2/q}{2}}\|f\|_{\mL_T^q\bC_\bba^\a},\no
\end{align}
Estimate \eqref{exies} follows by taking limits to the above inequality. For estimate \eqref{exies2} it follows similarly by (ii) of Lemma \ref{kryloves}.

Below we prove that $({\mathfrak F},Z,W)$ is a weak solution of SDE \eqref{3}. Note that by Proposition \ref{limit},
\begin{align}
\lim_{n\to\infty}\int_0^t b_{1;n}(r,Z_r)\dif r=:A_t^{b_1}\ \text{in}\ L^2.\no
\end{align}
It suffices to show that
\begin{align}\label{thm4.10:01}
\lim_{n\to\infty}\bE\Big|\int_0^t b_{1;n}(r,Z_r^n)\dif r-\int_0^t b_{1;n}(r,Z_r)\dif r\Big|^2=0,
\end{align}
and
\begin{align}\label{thm4.10:02}
\lim_{n\to\infty}\bE\Big|\int_0^t b_2(r,Z_r^n)\dif r-\int_0^t b_2(r,Z_r)\dif r\Big|^2=0.
\end{align}
First, we consider \eqref{thm4.10:01}. By Lebesgue's dominated convergence theorem, for fixed $m\in\mN$ we have 
$$\lim_{n\to\infty}{\bE}\Big|\int_0^t b_{1;m}(r,{Z_r^n})\dif r-\int_0^t b_{1;m}(r,{Z_r})\dif r\Big|^2=0.$$
Hence, to prove \eqref{thm4.10:01} we only need to show that
\begin{align}\label{thm4.10:03}
\lim_{m, m'\to\infty}\sup_{n}{\bE}\Big|\int_0^t (b_{1;m}-b_{1;m'})(r,{ Z_r^n})\dif r\Big|^2=0.
\end{align}
For $R>|z_0|$ and $n\in\mN$, define a stopping time 
$$\tau_R^n:=\inf\{t>0:|{Z_t^n}|\geq R\}.$$
Let $\chi\in C_c^\infty(\mR^{2d})$ with $\chi\equiv1$ on $B_1$ and $\chi\equiv0$ on $B_2^c$. Define a cut-off function
$$\chi_R(z):=\chi(z/R).$$
By Lemma \ref{kryloves}, we obtain
$$\lim_{m,m'\to\infty}\sup_n{\bE}\Big|\int_0^{t\wedge\tau_R^n} (b_{1;m}-b_{1;m'})(r,{Z_r^n})\dif r\Big|^2\lesssim\lim_{m,m'\to\infty}\|(b_{1;m}-b_{1;m'})\chi_R\|^2_{\mL_t^{q_b}\bC_{\bba}^{\a_b}}=0.$$
Additionally, by Chebyshev's inequality and \eqref{unimomentes}, we have
$$\lim_{R\to\infty}\sup_{n}{\bP}(\tau_R^n\leq T)=\lim_{R\to\infty}\sup_{n}{\bP}\Big(\sup_{t\in[0,T]}|{Z_t^n}|\geq R\Big)\leq\lim_{R\to\infty}\sup_{n}{\bE}\Big(\sup_{t\in[0,T]}|Z_t^n|\Big)/R=0.$$
Combining these limits, we obtain \eqref{thm4.10:03}, and therefore \eqref{thm4.10:01}. Next for \eqref{thm4.10:02}, denote $b_2^R(t,z):=b_2(t,z){\bf1}_{\{|z|<R\}}$. Similarly, we have
\begin{align}
\lim_{n\to\infty}\bE\Big|\int_0^t b_2^R(r,Z_r^n)\dif r-\int_0^t b_2^R(r,Z_r)\dif r\Big|^2=0,\no
\end{align}
which implies \eqref{thm4.10:02} by \eqref{unimomentes} as above.
\epf

\subsection{Uniqueness of  subcritical kinetic SDEs}\label{sec43}

The following lemma will play a crucial role in the proof of uniqueness.

\bl\label{geito}
Suppose \textbf{\emph{(H$^\text{sub}$)}} holds with $b_2=0$. Let $T>0$, $\a\in(-1,\a_b]$, $q\in(\ff{2}{1+\a},q_b]$ and $f\in\mL_T^q\bC_{\bba}^{\a}$. Let $u$ be the unique weak solution of the following backward PDE:
$$\partial_t u+\Delta_v u+v\cdot\nn_x u+b_1\odot\nabla_v u-\div_v b_1\preceq u=f,\ u(T)=0,$$
and let $(\mathfrak F,Z,W)$ be a weak solution of SDE \eqref{1} in the sense of Definition \ref{weaksol}. Then, for all $t\in[0,T]$, 
$$u(t,Z_t)=u(0,Z_0)+A_t^f+\sqrt{2}\int_0^t \nn_v u(s,Z_s)\dif W_s,$$
where $ A_t^f$ is defined in Proposition \ref{limit}.
\el

\bpf
By Theorem \ref{pdees}, there is a unique weak solution $u$ to the above backward PDE, with the regularity that for all $\t\in[0,2-\ff{2}{q})$,
\begin{align}\label{control}
(1+\l)^{\ff{\t}{2}}\|u\|_{\mL_T^\infty\bC_{\bba}^{2+\a-2/q-\t}}\leq C\|f\|_{\mL_T^q\bC_{\bba}^\a}.
\end{align}
Letting $u_n(t,z):=u(t,\cdot)*\G_n(z)$, by Remark \ref{remark3.3} we have 
\begin{align}\label{pdeap}
\partial_t u_n+\Delta_v u_n+v\cdot\nn_x u_n=f_n+g_n+h_n,
\end{align}
where
\begin{align}f_n:=f*\G_n,\ g_n:=v\cdot\nn_x u_n-(v\cdot\nn_x u)*\Gamma_n,\no
\end{align}
and
\begin{align}
h_n=-(b_1\odot\nabla_v u-\div_v b_1\preceq u)*\Gamma_n.\no
\end{align}
Taking $g=\nn_v u_n$ in \eqref{4.3}, by the generalized It\^{o}'s formula (cf. \cite{Follmer81}, \cite[Lemma 3.6]{ZZ17}) and \eqref{pdeap}, we obtain
\begin{align}\label{Ito}
u_n(t,Z_t)=u_n(0,Z_0)+A_t^{b_1\cdot\nn_v u_n+h_n}+\int_0^t(f_n+g_n)(s,Z_s)\dif s+\sqrt{2}\int_0^t \nn_v u_n(s,Z_s)\dif W_s.
\end{align}
Next, we estimate the difference between $u_n$ and $u$. For $\t'\in(0,2+\a-\ff{2}{q})$, 
\begin{align}\begin{split}
\|u_n-u\|_{\mL_T^\infty}&\leq\sup_{t\in[0,T]}\sup_{x,v}\int_{\mR^{2d}}|u(t,x-y,v-w)-u(t,x,v)|\G_n(y,w)\dif y\dif w\\\nonumber
&\leq\|u\|_{\mL_T^\infty\bC_{\bba}^{\t'}}\int_{\mR^{2d}}(|y|^{\ff{1}{3}}+|w|)^{\t'}\G_n(y,w)\dif y\dif w\overset{\eqref{control}}{\lesssim}n^{-\ff{\t'}{3}}.\nonumber
\end{split}\end{align}
Thus,
\begin{align}\label{es1}
(u_n(t,Z_t),u_n(0,Z_0))\rightarrow(u(t,Z_t),u(0,Z_0))\ \ a.s.,\ n\rightarrow\infty.
\end{align}
Similarly, there exists $\t''\in(0,1+\a-\ff{2}{q})$ such that
\begin{align}
\|\nn_v u_n-\nn_v u\|_{\mL_T^\infty}\lesssim n^{-\ff{\t''}{3}}\|u\|_{\mL_T^\infty\bC_\bba^{1+\t''}}\overset{\eqref{control}}{\lesssim}n^{-\ff{\t''}{3}}.\no
\end{align}
By BDG's inequality, we have
\begin{align}\label{es6}
\lim_{n\to\infty}\bE\Big|\int_0^t (\nn_v u_n(s,Z_s)-\nn_v u(s,Z_s))\dif W_s\Big|^2\lesssim\lim_{n\to\infty}\bE\Big(\int_0^t|(\nn_v u_n-\nn_v u)(s,Z_s)|^2\dif s\Big)=0.
\end{align}
Similarly, we note that
$$g_n(t,z)=\left(v\cdot\nn_x u_n-(v\cdot\nn_x u)*\G_n\right)(t,z)=\int_{\mR^{2d}}u(t,x-y,v-w)w\cdot\nn_y\G_n(y,w)\dif y\dif w,$$
which implies that
$$\|g_n(t)\|_\infty\lesssim n^{-2}\|u(t)\|_\infty\lesssim n^{-2}.$$
Hence,
\begin{align}\label{es2}
\lim_{n\to\infty}\int_0^t |g_n|(s,Z_s)\dif s=0.
\end{align}
In addition, by Proposition \ref{limit} and Theorem \ref{existence}, we have
\begin{align}\label{es3}
\lim_{n\to\infty}\bE\Big|\int_0^t f_n(s,Z_s)\dif s-A_t^f\Big|^2=0.
\end{align}
Finally, we decompose
$$b_1\cdot\nn_v u_n+h_n=h_1^n+h_2^n,$$
where
$$h_1^n:=b_1\odot\nn_v u_n-\G_n*(b_1\odot\nn_v u),$$
$$h_2^n:=\div_v b_1\preceq u_n-\G_n*(\div_v b_1\preceq u).$$
By Krylov's estimate \eqref{exies}, we have 
\begin{align}
\bE|A_t^{h_1^n}|\lesssim\|h_1^n\|_{\mL_T^{q_b}\bC_{\bba}^{\a_b}}\lesssim\|b_1\odot\nn_v(u_n-u)\|_{\mL_T^{q_b}\bC_{\bba}^{\a_b}}+\|b_1\odot\nn_v u-\G_n*(b_1\odot\nn_v u)\|_{\mL_T^{q_b}\bC_{\bba}^{\a_b}}.\nonumber
\end{align}
Note that there is an $\varepsilon>0$ such that for all $n\in\mN$,
$$\|b_1\odot\nn_v(u_n-u)\|_{\mL_T^{q_b}\bC_{\bba}^{\a_b}}\overset{\eqref{bonycontrol1}}{\lesssim}\|b_1\|_{\mL_T^{q_b}\bC_{\bba}^{\a_b}}\|u_n-u\|_{\mL_T^\infty\bC_{\bba}^{1}}\lesssim n^{-\varepsilon}.$$
Hence, by property of convolution we have
\begin{align}\label{es4}
\lim_{n\to\infty}\bE|A_t^{h_1^n}|=0.
\end{align}
Similarly, we can apply Krylov's estimate to $h_2^n$, yielding 
{\small
\begin{align}
\bE|A_t^{h_2^n}|\lesssim\|h_2^n\|_{\mL_T^{q_b}\bC_{\bba}^{\a_b}}\lesssim\|\div_v b_1\preceq u-\G_n*(\div_v b_1\preceq u)\|_{\mL_T^{q_b}\bC_{\bba}^{\a_b}}+\|\div_v b_1\preceq (u-u_n)\|_{\mL_T^{q_b}\bC_{\bba}^{\a_b}}.\nonumber
\end{align}}
Note that 
\begin{align}\begin{split}
\lim_{n\to\infty}\|\div_v b_1\preceq (u-u_n)\|_{\mL_T^{q_b}\bC_{\bba}^{\a_b}}&\lesssim\lim_{n\to\infty}\|\div_v b_1\preceq (u-u_n)\|_{\mL_T^{q_b}\bC_{\bba}^0}\\\nonumber
&\lesssim\|\div_v b_1\|_{\mL_T^{q_b}\bC_\bba^{\a_b}}\lim_{n\to\infty}\|u_n-u\|_{\mL_T^\infty\bC_{\bba}^{1}}=0.\nonumber
\end{split}\end{align}
As above we obtain
\begin{align}\label{es5}
\lim_{n\to\infty}\mE|A_t^{h_2^n}|=0.
\end{align}
Taking limits as $n\to\infty$ for both sides of \eqref{Ito} and by \eqref{es1}, \eqref{es6}, \eqref{es2}, \eqref{es3}, \eqref{es4}, \eqref{es5} and Proposition \ref{limit}, we complete the proof.
\epf

\bt\label{uniqueness}
Under \textbf{\emph{(H$^\text{sub}$)}} with $b_2=0$, the uniqueness of weak solution in Theorem \ref{existence} holds in the class that for any $\a\in(-1,\a_b]$, $q\in(\ff{2}{1+\a},q_b]$ and $p\geq2$, $(Z,\mL_T^q\bC_\bba^\a)$ satisfies the Krylov's estimate \eqref{krylov} with parameters $(p,1+\a-\ff{2}{q})$ and constant $C:=C(\Xi',\a,q,p)>0$.
\et

\bpf
Let $Z^1$ and $Z^2$ be two weak solutions of SDE \eqref{3} with the same initial distribution. By Lemma \ref{geito} with $T=t$, for any bounded continuous function $f$, we have
$$\bE\int_0^tf(Z_s^1)\dif s=-\bE u(0,Z_0^1)=-\bE u(0,Z_0^2)=\bE\int_0^tf(Z_s^2)\dif s.$$
Thus, $\bE f(Z_t^1)=\bE f(Z_t^2)$, which implies that $Z^1$ and $Z^2$ have the same one-dimensional marginal distribution. Hence, the weak uniqueness follows by a standard method (cf. \cite[Theorem 4.4.3]{EK86}), and we omit the details here.
\epf

We now extend the result to the general case where $b=b_1+b_2$. The main difficulty lies in the fact that the previous method is no longer applicable. To overcome this, we employ a standard localization technique as in \cite{SV06}. We begin by introducing the notion of a martingale solution.

\bd[Martingale Solutions]\label{mart sol}
For $(s,z)\in[0,T]\times\mR^{2d}$, a probability measure $\mP\in\sP(\mC_T)$ is called a martingale solution to SDE \eqref{3}, with drift $b=b_1+b_2$ where $b_1$ belongs to certain anisotropic H\"{o}lder space and $b_2$ is sufficiently regular, starting from $z$ at time $s$ if
\begin{itemize}
\item[(i)] $\mP(\om_t=z,\forall t\in[0,s])=1$.
\item[(ii)] For any $f\in C^2(\mR^{2d})$, 
\begin{align}
M_t^f:=f(\om_t)-f(\om_s)-\int_s^t(\Delta_v+v\cdot\nn_x+b_2(r)\cdot\nn_v)f(\om_r)\dif r-A_{s,t}^{b_1\cdot\nn_v f},\ s\leq t\leq T\no
\end{align}
is a $\sB_t$-martingale under $\mP$, where $A_{s,t}^{b_1\cdot\nn_v f}:=\lim_{n\to\infty}\int_s^t (b_{1;n}(r)\cdot\nn_v f)(\om_r)\dif r$ exists in $L^2$-sense.
\end{itemize}
The set of all the above martingale solutions starting from $z$ at time $s$ is denoted by $\cM_z^s$.
\ed

\br
By an argument similar to that of \cite[Proposition 3.13]{ZZ17}, for any $(s,z_0)\in[0,T]\times\mR^{2d}$, we have $\mP\in\cM_{z_0}^s$ if and only if there exists a weak solution $(\mathfrak F,Z,W)$ to SDE \eqref{3} starting from $z_0$, with the drift term $b(t,z)$ replaced by $b(t+s,z)$, in the sense of Definition \ref{weaksol} such that ${\bf P}\circ Z^{-1}=\mP\circ\t_s^{-1}$, where $\t_s$ is the standard shift operator defined by $\t_s(\om_t)=\om_{t+s}$.
\er

A standard truncated method will be applied to derive the uniqueness. Let $\wt\cM_z^s(R)$ denote the set of martingale solutions of SDE \eqref{3} obtained by only replacing $b_2$ by $b_2^R:=b_2{\bf1}_{\{|z|<R\}}$ in Definition \ref{mart sol}. By Theorem \ref{uniqueness}, we have

\bl\label{lemma4.11}
Under \textbf{\emph{(H$^\text{sub}$)}}, for each $R\in(0,\infty)$ and $(s,z)\in[0,T]\times\mR^{2d}$, there exists a unique element in $\wt\cM_z^s(R)$.
\el

The following two lemmas are from \cite[Lemma 6.1.1 and Theorem 6.1.2]{SV06}.

\bl
Fix $s>0$ and $\eta\in\mC_T$. Suppose that $\mQ\in\sP(\mC_T)$ satisfies $\mQ(\{\omega\in\mC_T:\omega(s)=\eta(s)\})=1$. Then there exists a unique probability measure $\delta_\eta\otimes_s\mQ$ on $\mC_T$ such that
\begin{align}
\delta_\eta\otimes_s\mQ(\{\omega\in\mC_T:\omega(t)=\eta(t),\forall t\in[0,s]\})=1,\no
\end{align}
and
\begin{align}
\delta_\eta\otimes_s\mQ=\mQ\ \text{on}\ \sB^s.\no
\end{align}
\el

\bl\label{svlemma}
Let $\tau$ be a finite $\sB_t$-stopping time on $\mC_T$. Suppose that $\mQ:\eta\mapsto \mQ_\eta$ is a map from $\mC_T$ to $\sP(\mC_T)$ such that for each $A\in\sB_T$, $\eta\mapsto \mQ_\eta(A)$ is $\sB_\tau$-measurable and for any $\eta\in\mC_T$,
\begin{align}
\mQ_\eta(\{\om\in\mC_T:\om(\tau(\eta))=\eta(\tau(\eta))\})=1.\no
\end{align}
Let $\mP\in\sP(\mC_T)$. Then
\begin{itemize}
\item[(i)] There exists a unique $\mP\otimes_\tau \mQ\in\sP(\mC_T)$ such that $\mP\otimes_\tau \mQ=\mP$ on $\sB_\tau$ and $\{\de_\eta\otimes_{\tau(\eta)}\mQ_\eta\}_{\eta\in\mC_T}$ is a regular conditional probability distribution of $\mP\otimes_\tau \mQ$ given $\sB_\tau$.
\item[(ii)] If $M:\mR_+\times\mC_T\to\mR$ is progressive measurable and right-continuous such that $M_t^\tau:=M_{t\wedge\tau}$ is a $\mP$-martingale and $M_t-M_t^{\tau(\eta)}$ is a $\mQ_\eta$-martingale for each $\eta\in\mC_T$, then $M_t$ is a $\mP\otimes_\tau \mQ$-martingale.
\end{itemize}
\el

The following lemma can be obtained by a standard stopping time argument together with Lemma \ref{svlemma} (cf. \cite[Lemma 5.9]{CZZ21}).

\bl\label{lemma4.13}
Fix $R\in(0,\infty)$ and $(s,z)\in\mR_+\times\mR^{2d}$. Define the $\sB_t$-stopping time
$$\tau(\om):=\inf\{t>s:|\om(t)|>R\},$$
and for $f\in C^2(\mR^{2d})$,
\begin{align}
\wt M_t^f:=f(\om_t)-f(\om_s)-\int_s^t(\Delta_v+v\cdot\nn_x+b_2^R(r)\cdot\nn_v)f(\om_r)\dif r-A_{s,t}^{b_1\cdot\nn_v f},\ t\geq s.\no
\end{align}
Suppose that $\mP\in\cM_z^s$. Then $(\wt M^f)_t^\tau:=\wt M_{t\wedge\tau}^f$ is a $\sB_t$-martingale under $\mP$. Moreover, if we denote $\mQ_\eta:=\wt \mP_{\tau(\eta),\eta_{\tau(\eta)}}$ for $\eta\in\mC_T$, where $\wt \mP_{s,z}$ is the unique element in $\wt\cM_z^s(R)$, then $\mP\otimes_\tau\mQ\in\wt\cM^s_z(R)$.
\el

We prove the following theorem by arguments similar to those in \cite[Theorem 10.4]{Pinsky95}.

\bt\label{uniuni}
Let $\tau_R=\inf\{t>s;|\om(t)|>R\}$ (with $\inf\emptyset=\infty$ by convention). If $\lim_{R\to\infty}\tau_R=\infty$, then the uniqueness of $\wt\cM_z^s(R)$ for all $R\in(0,\infty)$ and $(s,z)\in\mR_+\times\mR^{2d}$ implies the uniqueness of $\cM_z^s$ for all $(s,z)\in\mR_+\times\mR^{2d}$.
\et

\bpf
Assume $\mP^{(1)},\mP^{(2)}\in\cM_z^s$. Denote $\mQ:\eta\mapsto\mQ_\eta:=\wt \mP_{\tau_R(\eta),\eta_{\tau_R(\eta)}}$ for $\eta\in\mC_T$, where $\wt \mP_{s_0,z_0}$ is the unique element in $\wt\cM_{z_0}^{s_0}(R)$. By Lemma \ref{lemma4.13} we have
\begin{align}
\mP^{(i)}\otimes_{\tau_R} \mQ\in\wt\cM^s_z(R),\ i=1,2.\no
\end{align}
By the uniqueness of $\wt\cM_z^s(R)$, we have $\mP^{(i)}\otimes_{\tau_R} \mQ=\wt\mP_{s,z}$ for $i=1,2$. It then follows from Lemma \ref{svlemma} (i) that
\begin{align}
\mP^{(1)}=\mP^{(1)}\otimes_{\tau_R}\mQ=\wt\mP_{s,z}=\mP^{(2)}\otimes_{\tau_R} \mQ=\mP^{(2)}\ \text{on}\ \sB_{\tau_R}.\no
\end{align}
Thus, letting $R\to\infty$ we have
\begin{align}
\mP^{(1)}=\mP^{(2)},\no
\end{align}
which completes the proof.
\epf

Now we can show the well-posedness of SDE \eqref{3}.

\bt\label{wellpodesness}
Under \textbf{\emph{(H$^\text{sub}$)}}, there exists a unique weak solution $(\mathfrak F,Z,W)$ to SDE \eqref{3} starting from $z_0$ such that the following Krylov's estimate holds: for any $f\in\mL_T^1 C_b\cap\mL_T^q\bC_\bba^\a$ with $\a\in(-1,\a_b]$, $q\in(\ff{2}{1+\a},q_b]$, $p\geq2$, and stopping times $0\leq \tau_0\leq \tau_1\leq T$, $\tau_1-\tau_0\leq\s$,
\begin{align}
\Big\|\int_{\tau_0}^{\tau_1}f(s,Z_s)\dif s\Big\|_{L^p(\Omega)}\lesssim_C \s^{\ff{2+\a-2/q}{2}}\|f\|_{\mL_T^q\bC_\bba^\a},\no
\end{align}
where the constant $C$ is independent of $t_0$, $t_1$ and $f$.  Additionally, for any $\de\in\mR$,
\begin{align}
\bE\bigg(\sup_{t\in[0,T]}\rho_\de(Z_t)\bigg)\lesssim_{C(\de)}\rho_\de(z_0).\no
\end{align}
\et

\bpf
The existence of solutions follows immediately from Theorem \ref{existence}. Regarding uniqueness, Lemma \ref{lemma4.11} establishes the uniqueness of $\wt\cM_z^s(R)$ for all $R\in(0,\infty)$ and $(s,z)\in\mR_+\times\mR^{2d}$. Furthermore, \eqref{unimomentes} implies that $\lim_{R\to\infty}\tau_R=\infty$. Applying Theorem \ref{uniuni}, we consequently obtain the uniqueness of $\cM_z^s$ for all $(s,z)\in\mR_+\times\mR^{2d}$. The proof is complete by observing that ${\bf P}\circ Z^{-1}\in\cM_{z_0}^0$.
\epf

Finally we can give
\bpf[Proof of Theorem \ref{sec1:main2}]
By Remark \ref{GubinelliDecomposition}, for any $b\in\mL_T^{q_b}\bC_\bba^{\a_b}$ with $\k\in[0,1+\a_b)$, assumption \textbf{(H$^\text{sub}$)} is satisfied. In this case, the well-posedness of the SDE and estimate \eqref{sec1:main2:es1} follow directly from Theorem~\ref{wellpodesness}. Furthermore, when $\k\in[0,\ff{1+\a_b-2/q_b}{3+\a_b-2/q_b}]$, assumption \textbf{(H$_\text{w}^\text{sub}$)} holds, and estimate \eqref{sec1:main2:es2} follows immediately from \eqref{exies2}.
\epf

\section{Applications to Gaussian fields}\label{sec5}
In this section, we construct a concrete example of a Gaussian noise $b$ that satisfies the assumptions of Theorem \ref{sec1:main2}. Our construction proceeds in two steps: first, we establish the regularity properties stated in Theorem \ref{GFFre}, and subsequently, we explicitly construct the desired noise.

We begin by introducing the Gaussian noise with spectral measure $\mu$. Let $\mu$ be a signed Radon measure on $\mR^{2d}$ satisfying the following conditions:

\begin{enumerate}[\textbf{(A)}]

\item $\mu$ is a symmetric measure, meaning that
\begin{align}
\mu(\dif\xi, \dif\eta)=\mu(-\dif\xi, \dif\eta)=\mu(\dif\xi, -\dif\eta).\no
\end{align}
Moreover, for some $\ell\in\mR$, it holds that
\begin{align}
\int_{\mR^{2d}} (1+|\zeta|_\bba)^\ell |\mu|_{\text{var}}(\dif\zeta)<\infty,\no
\end{align}
where $|\cdot|_{\text{var}}$ denotes the total variation measure of a signed measure.
\end{enumerate}

Let $\mH$ be the completion of $\sS(\mR^{2d};\mR)$ with respect to the inner product
\begin{align}
\langle f,g\rangle_\mH:=\int_{\mR^{2d}}\hat f(\zeta)\overline{\hat g(\zeta)}\mu(\dif\zeta)<\infty.\no
\end{align}
Now we introduce the Gaussian random field on $\mH$ (see \cite{Jason97}).

\bd
Let $U$ be the real-valued Gaussian random field on $\mH$; that is, $U$ is a continuous linear operator from $\mH$ to $L^2(\Omega,\bP)$, and for any $f\in\mH$, $U(f)$ is a real-valued Gaussian random variable with mean zero and variance $\Vert f\Vert_\mH$. Moreover, 
\begin{align}\label{Gaucov}
\bE(U(f)U(g))=\langle f,g\rangle_\mH.
\end{align}
We call $U$ the Gaussian noise with spectral measure $\mu$.
\ed

Below, for a Banach space $\mB$, denote by $\bC_\bba^\a(\mB)$ the $\mB$-valued H\"{o}lder space. Recall $\rho_\k(z)=((1+|x|^2)^{1/3}+1+|v|^2)^{-\k/2}$ for $\k\in\mR$. Let $\bC_\bba^\a(\rho_\k;\mB)$ be the weighted $\mB$-valued H\"{o}lder space with norm (see \cite[Theorem 2.7]{HZZZ24})
\begin{align}
\|f\|_{\bC_\bba^\a(\rho_\k;\mB)}:=\|f\rho_\k\|_{\bC_\bba^\a(\mB)}.\no
\end{align}
The following main result addresses the regularity of $U$. 

\bt\label{GFFre}
Under \emph{\textbf{(A)}}, for any $p\in[1,\infty)$, $\k>\ff{3d}{p}$ and $\b<\ff{\ell}{2}-\ff{4d}{p}$, it holds that
\begin{align}
U\in L^p(\Omega;\bC_\bba^\b(\rho_\k;\mR)).
\end{align}
\et

The key point of Theorem \ref{GFFre} is to view $U$ as a random distribution of the variable $z:=(x,v)\in\mR^{2d}$. To formalize this, we introduce the following notation. For any given $\varphi\in\mH$, we define
\begin{align}
(\varphi*U)(z):=U_\varphi(z):=U(\tau_z\varphi),\no
\end{align}
where $\tau_z\varphi(\cdot):=\varphi(\cdot-z)\in\mH$. By \cite[Lemma B.1]{HZZZ24}, the mapping $z\mapsto U_\varphi(z)$ has a smooth version.

\bpf[Proof of Theorem \ref{GFFre}]
By \eqref{supp}, for all $\g\in\mR$, it follows that 
\begin{align}
|\phi_j^\bba(\zeta)|\lesssim 1\wedge(2^{\g j}(1+|\zeta|_\bba)^{-\g}),\ j\geq0,\ \zeta\in\mR^{2d}.\no
\end{align}
By the hypercontractivity of Gaussian random variables, for any fixed $z\in\mR^{2d}$, $p\in[1,\infty)$ and $j\in\mN$, we have 
\begin{align}
&\|\cR_j^\bba U(z)\|_{L^p(\Omega)}\lesssim\|\cR_j^\bba U(z)\|_{L^2(\Omega)}=\|U(\tau_z\check\phi_j^\bba)\|_{L^2(\Omega)}\no\\
&\qquad\overset{\eqref{Gaucov}}{=}\Big(\int_{\mR^{2d}}|\phi_j^\bba(\zeta)|^2\mu(\dif\zeta)\Big)^{1/2}\lesssim2^{-\ff{\ell j}{2}}\Big(\int_{\mR^{2d}}(1+|\zeta|_\bba)^\ell\mu(\dif \zeta)\Big)^{1/2}\lesssim2^{-\ff{\ell j}{2}}.\no
\end{align}
For $\k>\ff{3d}{p}$,
\begin{align}
\|\rho_\k\cR^\bba_j U\|^p_{L^p(\mR^{2d};L^p(\Omega))}=\int_{\mR^{2d}}\rho_{p\k}(z)\|\cR_j^\bba U(z)\|^p_{L^p(\Omega)}\dif z\lesssim2^{-\ff{p\ell j}{2}}.\no
\end{align}
Then for $\b<\ff{\ell}{2}-\ff{4d}{p}$, by Bernstein's inequality \eqref{bernstein2} and Fubini's Theorem we have
\begin{align}
&\Big\|\sup_{j\in\mN}2^{j\b}\|\rho_\k\cR_j^\bba U\|_{L^\infty(\mR^{2d})}\Big\|_{L^p(\Omega)}\lesssim\Big\|\sup_{j\in\mN}2^{j(\b+\ff{4d}{p})}\|\rho_\k\cR_j^\bba U\|_{L^p(\mR^{2d})}\Big\|_{L^p(\Omega)}\no\\
&\qquad\leq\sum_{j\in\mN}2^{j(\b+\ff{4d}{p})}\|\rho_\k\cR_j^\bba U\|_{L^p(\mR^{2d};L^p(\Omega))}\lesssim\sum_{j\in\mN}2^{j(\b+\ff{4d}{p}-\ff{\ell}{2})}<\infty.\no
\end{align}
The proof is complete.
\epf

Now we provide the example to illustrate our results.

\bx\label{eg5.3}
For $\g\in[0,d)$, let 
\begin{align}
\mu(\dif\xi,\dif\eta)=|\xi|^{-\g}\dif\xi\delta_0(\dif\eta),\no
\end{align}
where $\dif\xi$ denotes the common Lebesgue measure on $\mR^{d}$ and $\delta_0(\dif\eta)$ is the Dirac measure concentrated at 0. For any real-valued Gaussian noise $U$ with above spectral noise $\mu$, by \cite[p. 117, Lemma 2]{Stein70}, there exists some constant $c_{d,\g}>0$ such that
\begin{align}
\hat\mu(z)=\hat\mu(x,v)=c_{d,\g}|x|^{\g-d},\ z=(x,v).\no
\end{align}
In particular, for any $f,g\in\sS(\mR^{2d})$, by \eqref{Gaucov} we have
\begin{align}
\mE(U(f)U(g))&=\int_{\mR^{2d}}\hat f(\zeta)\hat g(-\zeta)\mu(\dif\zeta)=\int_{\mR^{2d}}\int_{\mR^{2d}}f(z)g(z')\hat\mu(z-z')\dif z\dif z'\no\\
&=c_{d,\g}\int_{\mR^d}\int_{\mR^d}\Big(\int_{\mR^d}f(x,v)\dif v\Big)\Big(\int_{\mR^d}g(x',v')\dif v'\Big)\ff{\dif x\dif x'}{|x-x'|^{d-\g}}.\no
\end{align}
Generally, smaller $\g$ implies worse regularity of Gaussian random field $U$. Moreover, by the change of variables we obtain
\begin{align}
\|\cR_j^\bba U(z)\|_{L^p(\Omega)}\lesssim\Big(\int_{\mR^{2d}}|\phi_j^\bba(\zeta)|^2\mu(\dif\zeta)\Big)^{1/2}=2^{\ff{3j(d-\g)}{2}}\Big(\int_{\mR^d}|\phi_0^\bba(\xi,0)|^2\ff{\dif\xi}{|\xi|^\g}\Big)^{1/2}.\no
\end{align}
Thus, for $\k>\ff{3d}{p}$ and $\b<\ff{3(\g-d)}{2}-\ff{4d}{p}$ we conclude that
$$U\in L^p(\Omega;\bC_\bba^\beta(\rho_\k;\mR)).$$
Next we construct a Gaussian field that is divergence free in  $v$. Fix $\varphi\in\sS(\mR^d)$ such that $\int_{\mR^d}\varphi=1$. For any $f\in\sS(\mR^d)$, we define
\begin{align}
\wt U(f):=U(\wt f),\ \wt f(x,v):=f(x)\varphi(v).\no
\end{align}
In this case, for any $f,g\in\sS(\mR^d)$,
\begin{align}
\mE(\wt U(f)\wt U(g))=c_{d,\g}\int_{\mR^d}\int_{\mR^d}f(x)g(x')\ff{\dif x\dif x'}{|x-x'|^{d-\g}}=\int_{\mR^d}\hat f(\xi)\hat g(\xi)\ff{\dif\xi}{|\xi|^\g}.\no
\end{align}
Let ${\bf U}=(\wt U_1,\cdots,\wt U_d)$ be d-independent Gaussian noises defined above. Then it holds that for any $\k>\ff{3d}{p}$ and $\b<\ff{3(\g-d)}{2}-\ff{4d}{p}$,
$${\bf U}\in L^p(\Omega;\bC_\bba^\beta(\rho_\k;\mR^d)).$$
In particular, one can choose $p$ big enough such that $\k>\ff{3d}{p}$ and $\b+\ff{4d}{p}<\ff{3(\g-d)}{2}$, and for $\bP$-almost all $\omega$,
$${\bf U}(\omega,\cdot)\in\bC_\bba^\b(\rho_\k),$$
and
$$\div_v {\bf U}(\omega,\cdot)=0.$$
Hence, for any $\g\in(d-\ff{2}{3},d)$, it follows from Theorem \ref{sec1:main2} that there exists a unique weak solution to the following kinetic SDE:
\begin{align}
\dif X_t=V_t\dif t,\ \dif V_t={\bf U}(\omega,\cdot)(X_t)\dif t+\sqrt{2}\dif W_t.\no
\end{align}
\ex

\end{document}